\documentclass[twoside,10pt,a4paper,leqno]{article}
\usepackage{amssymb,amsmath,amscd,euscript,verbatim,array}
\usepackage[center,pagestyles]{titlesec}
\usepackage{geometry}
\usepackage{anysize}
\usepackage{fancyhdr}
\usepackage{indentfirst}
\usepackage{graphicx}
\usepackage{color}
\usepackage{ifpdf}
\usepackage{mathrsfs}
\usepackage{cite}
\usepackage[numbers,sort&compress]{natbib}
\marginsize{3.5cm}{3.5cm}{2.75cm}{2.75cm}

\titleformat{\section}{\centering\normalsize}{\thesection.}{0.5em}{}
\titleformat{\subsection}{\normalsize\bfseries}{\thesubsection.}{0.5em}{}
\titleformat{\subsubsection}{\normalsize\bfseries}{\thesubsubsection.}{0.5em}{}
\newcommand{\N}{\mathbb{N}}

\newcommand{\R}{\mathbb{R}}

\newtheorem{Theorem}{Theorem}[section]

\newtheorem{Lemma}[Theorem]{Lemma}
\newtheorem{Exercise}[Theorem]{Exercise}
\newtheorem{Proposition}[Theorem]{Proposition}
\newtheorem{Remark}[Theorem]{Remark}

\newcommand{\gm}{\gamma}

\newcommand{\bthm}{\begin{Theorem}}
\newcommand{\ethm}{\end{Theorem}}
\newcommand{\bpr}{\begin{Proposition}}
\newcommand{\epr}{\end{Proposition}}
\newcommand{\blm}{\begin{Lemma}}
\newcommand{\elm}{\end{Lemma}}
\newcommand{\bex}{\begin{Exercise}}
\newcommand{\eex}{\end{Exercise}}
\newcommand{\be}{\begin{equation}}
\newcommand{\ee}{\end{equation}}
\newcommand{\beal}{\begin{aligned}}
\newcommand{\enal}{\end{aligned}}
\newcommand{\brm}{\begin{Remark}}
\newcommand{\erm}{\end{Remark}}
\newcounter{item}[section]

\newcommand{\Proof}{\textbf{Proof}\hspace{0.3cm}}
\newcommand{\End}{\ensuremath{\hfill{\Box}}\\}
\renewcommand{\title}[1]{\begin{center}\bf\large #1\end{center}}
\renewcommand{\author}[1]{\begin{center}\normalsize #1\end{center}}
\renewcommand{\date}[1]{\begin{center}#1\end{center}}

\makeatletter \@addtoreset{equation}{section}
\makeatother
\begin{document}
\vspace{10pt}
\title{VARIATIONAL PRINCIPLE FOR CONTACT TONELLI HAMILTONIAN SYSTEMS}

\vspace{6pt}
\author{\sc Lin Wang \& Jun Yan}

\vspace{10pt} \thispagestyle{plain}
%{\begingroup\makeatletter
%\let\@makefnmark\relax  \footnotetext{$\ast$ Corresponding author}
%\makeatother\endgroup}

\begin{quote}
\small {\sc Abstract.} We establish an implicit variational principle for the equations of the contact flow generated by the  Hamiltonian $H(x,u,p)$ with respect to the contact 1-form $\alpha=du-pdx$ under Tonelli and Osgood growth assumptions. It is the first step to  generalize Mather's global variational method from the Hamiltonian dynamics to the contact Hamiltonian dynamics.
\end{quote}
\begin{quote}
\end{quote}
\begin{quote}
\end{quote} \vspace{25pt}
\tableofcontents
\newpage
\section{\sc Introduction and main results}
Let $M$ be an $n$-dimensional closed manifold. The cotangent bundle $T^*M$ has a natural symplectic structure $\omega$. The pair $(T^*M,\omega)$ is  a symplectic manifold.   Any symplectic form is locally diffeomorphic to
\[dx\wedge dp=dx_1\wedge dp_1+\cdots+dx_n\wedge dp_n.\]
Let $H:T^*M\rightarrow\R$ be a $C^r$ ($r\geq 2$)  function called a Hamiltonian. The equations of the  Hamiltonian flow generated by $H(x,p)$ with respect to $\omega$ has the form in local coordinates:
\begin{equation}\label{hjech22}
\begin{cases}
\dot{x}=\frac{\partial H}{\partial p},\\
\dot{p}=-\frac{\partial H}{\partial x}.
\end{cases}
\end{equation}
In the early 1990's, based on the variational principle in Tonelli Lagrangian systems, J. Mather \cite{M1,Mvc}  founded a seminar work (so called Mather theory) on global action minimizing orbits generated by (\ref{hjech22}). Mather theory has a profound impact on many fields, such as Hamiltonian dynamics (Arnold diffusion \cite{Be2,Chengcq,CY1,CY2,KAZ}, Ma\~{n}\'{e} conjecture \cite{CFR,FR1,FR2}, converse KAM theory (\cite{B2,CC1,CW,Fo,M4,W1})), Hamilton-Jacobi equations (weak KAM theory \cite{Be,CI,F1,F2,F3,F22,FS,KA}) and symplectic topology (symplectic invariants \cite{Be1,PPS,SV})  {\it etc.}.

Comparatively, let $J^1(M,\R)$ denote the manifold of 1-jets of functions on $M$. If $(x_1,\ldots,x_n)$ are local coordinates on $M$, then $J^1(M,\R)$ is defined by a collection
\[(x_1,\ldots,x_n,u,p_1,\ldots,p_n),\]
where $p_i=(\partial u/\partial x_i)(x)$.
The standard contact form on $J^1(M,\R)$ is the 1-form
\[\alpha=du-pdx\quad(pdx=p_1dx_1+\cdots+p_ndx_n).\]
 $J^1(M,\R)$ has a natural contact structure, which is globally defined by the Pfaffian equation $\alpha=0$. The pair $(J^1(M,\R),\alpha)$ is  a contact manifold. Similar to symplectic manifolds, any contact form is locally diffeomorphic to $\alpha$ (Darboux's normal form). There is a canonical diffeomorphism between $J^1(M,\R)$ and  $T^*M\times\R$. Thus, $(T^*M\times\R,\alpha)$ is also a contact manifold. Let $H:T^*M\times\R\rightarrow\R$ be a $C^r$ ($r\geq 2$) function. The equations of the contact flow generated by $H$ with respect to $\alpha$ have the following form in  local coordinates \cite{Ar1}:
\begin{equation}\label{hjech}
\begin{cases}
\dot{x}=\frac{\partial H}{\partial p},\\
\dot{p}=-\frac{\partial H}{\partial x}-\frac{\partial H}{\partial u}p,\\
\dot{u}=\frac{\partial H}{\partial p}p-H.
\end{cases}
\end{equation}

In this paper, we are devoted to establishing a variational principle for the $2n+1$th coupled system (\ref{hjech}) under Tonelli and Osgood growth assumptions. By the analogy of the celebrated Mather theory, the variational principle can be viewed as a stepping stone to make a further exploration on contact Hamiltonian systems, general Hamilon-Jacobi equations (depending $u$ explicitly) \cite{SWY}  and contact topology {\it etc.}.

Precisely speaking, we are concerned with a  $C^r$ ($r\geq 2$) Hamiltonian $H(x,u,p)$ satisfying the following conditions:
\begin{itemize}
\item [\textbf{(H1)}] \textbf{Positive Definiteness}: For every $(x,u)\in M\times\R$, the second partial derivative $\partial^2 H/\partial p^2 (x,u,p)$ is positive definite as a quadratic form;
%\item [\textbf{(H2)}] \textbf{Superlinearity in the Fibers}: For every compact set $I$, $H(x,u,p)$ ($u\in I$) is uniformly superlinear growth with respect  to $p$;
    \item [\textbf{(H2)}] \textbf{Superlinearity in the Fibers}: For every $(x,u)\in M\times\R$, $H(x,u,p)$ is  superlinear growth with respect  to $p$;
\item [\textbf{(H3)}] \textbf{Completeness of the Flow}: The flows of (\ref{hjech}) generated by $H(x,u,p)$ are complete;
\item [\textbf{(H4)}] \textbf{Osgood Growth}: For every compact set $K\subset T^*M$, there exists a continuous function $f_K(u)$ defined on $[0,+\infty)$
    with $f_K(u)>0$ and  $\int_0^{+\infty}\frac{1}{f_K(u)}du=+\infty$  such that for any $(x,p,u)\in K\times\R^+$,
    \[H(x,u,p)\geq \langle \frac{\partial H}{\partial p}(x,u,p),p\rangle-f_K(u).\]
\end{itemize}
%It is easy to see that $\int_0^\infty\frac{1}{f_K(u)}du$ is divergent if and only if the flow generated by $\dot{u}=f_K(u)$ is complete \cite{Ar2}.
Literately, (H1)-(H4) are called Tonelli conditions \cite{F3,M1}. (H4)  can be referred as Osgood condition \cite{os}.

We use $\mathcal{L}: T^*M\rightarrow TM$ to denote  the Legendre transformation. Let
$\bar{\mathcal{L}}:=(\mathcal{L}, Id)$, where $Id$ denotes the identity map from $\R$ to $\R$. Then $\bar{\mathcal{L}}$ denote a diffeomorphism from $T^*M\times\R$ to $TM\times\R$. By $\bar{\mathcal{L}}$,
the Lagrangian $L(x,u, \dot{x})$ associated to $H(x,u,p)$ can be denoted by
\[L(x,u, \dot{x}):=\sup_p\{\langle \dot{x},p\rangle-H(x,u,p)\}.\]
Let $\Psi_t$ denote the flows of (\ref{hjech}) generated by $H(x,u,p)$. Let the flow $\Phi_t:=\bar{\mathcal{L}}\circ\Psi_t\circ\bar{\mathcal{L}}^{-1}$ defined in $TM\times\R$. We call $\Phi_t$ the flow generated by $L(x,u,\dot{x})$.
 Based on
(H1)-(H4), we have:
\begin{itemize}
\item [\textbf{(L1)}] \textbf{Positive Definiteness}: For every $(x,u)\in M\times\R$, the second partial derivative $\partial^2 L/\partial {\dot{x}}^2 (x,u,\dot{x})$ is positive definite as a quadratic form;
%\item [\textbf{(L2)}] \textbf{Superlinearity in the Fibers}: For every compact set $I$, $L(x,u,\dot{x})$ ($u\in I$) is uniformly superlinear growth with respect  to $\dot{x}$;
    \item [\textbf{(L2)}] \textbf{Superlinearity in the Fibers}: For every $(x,u)\in M\times\R$, $L(x,u,\dot{x})$ is superlinear growth with respect  to $\dot{x}$;
\item [\textbf{(L3)}] \textbf{Completeness of the Flow}: The flow $\Phi_t$  is complete;
\item [\textbf{(L4)}] \textbf{Osgood Growth}: For every compact set $K\subset TM$, there exists a continuous function $f_K(u)$ defined on $[0,+\infty)$
    with $f_K(u)>0$ and   $\int_0^{+\infty}\frac{1}{f_K(u)}du=+\infty$   such that for  any $(x,\dot{x},u)\in K\times\R^+$,
    \[L(x,u,\dot{x})\leq f_K(u).\]
\end{itemize}
It is easy to see that (L4) is more general than the  monotonicity (non-increasing)  and the uniform Lipschitzity  of $L$ with respect to $u$, for which $f_K(u)$ is corresponding to a constant function and an affine function respectively.

%To avoid the ambiguity, we denote the solution of (\ref{hjech})  by $(X(t),U(t),P(t))$.
If a Hamiltonian $H(x,u,p)$ satisfies (H1)-(H4), then we obtain the following three theorems:
\begin{Theorem}\label{two}
For given $x_0\in M$, $u_0\in\R$ and $T>0$, there exists a  unique continuous function $h_{x_0,u_0}(x,t)$ defined on $M\times (0,T]$ satisfying
\begin{equation}\label{kkkk}
h_{x_0,u_0}(x,t)=u_0+\inf_{\substack{\gm(t)=x \\  \gm(0)=x_0} }\int_0^tL(\gm(\tau),h_{x_0,u_0}(\gm(\tau),\tau),\dot{\gm}(\tau))d\tau,
\end{equation}where the infimum is taken among the continuous and piecewise $C^1$ curves $\gm:[0,t]\rightarrow M$. Moreover, the infimum is attained at a $C^1$ curve denoted by $\bar{\gm}$. Let
\[X(t):=\bar{\gm}(t),\quad U(t):=h_{x_0,u_0}(\bar{\gm}(t),t),\quad P(t):=\bar{\mathcal{L}}^{-1}(\bar{\gm}(t),h_{x_0,u_0}(\bar{\gm}(t),t),\dot{\bar{\gm}}(t)),\]
then $(X(t),U(t),P(t))$ satisfies the equation (\ref{hjech}).
\end{Theorem}

\begin{Theorem}\label{two2}
Let $\mathcal{S}_{x_0,u_0}^x(t)$ denote the set of the solutions of (\ref{hjech}) satisfying $X(0)=x_0$, $X(t)=x$ and $U(0)=u_0$, then we have
\begin{equation}
h_{x_0,u_0}(x,t)=\inf\left\{U(t):(X(s),U(s),P(s))\in \mathcal{S}_{x_0,u_0}^x(t)\right\}.
\end{equation}
\end{Theorem}

\begin{Theorem}\label{two1}
For each solution $(X(t),U(t),P(t))$ of  (\ref{hjech}), there exists $\epsilon>0$ such that for any $t\in (0,\epsilon]$,
\[U(t)=h_{x_0,u_0}(X(t),t),\]
where $X(0)=x_0$ and $U(0)=u_0$.
\end{Theorem}

Theorem \ref{two2} can be obtained following from the proof of Theorem \ref{two}. Under the assumptions (H1)-(H4), Theorem \ref{two} establishes a variational principle for the contact Hamiltonian equation (\ref{hjech});
Theorem \ref{two2} provides a representation of the minimal action $h_{x_0,u_0}(x,t)$ by certain global solutions of (\ref{hjech}) with initial  and boundary conditions; Theorem \ref{two1} shows that each short solution of (\ref{hjech}) is a strict minimizer in the sense of the variation. Theorem \ref{two1} generalizes an ancient result due to Weierstrass that sufficiently short solutions of the Euler-Lagrange equation are strict minimizers \cite{M1}.

%Traditionally, the main tool to address the contact topology is symplectification, so that some techniques from the symplectic topology can be used. However, some crucial properties, for instance the convexity of the contact Hamiltonian, could be lost due to symplectification. Based on the implicit variational principle, it is hopeful to handle the contact topology directly by dynamical approaches.

The paper is outlined as follows. In Section 2, in order to avoid the singularity of $h_{x_0,u_0}(x,t)$ at $t=0$ for $x_0\neq x$, the Lagrangian will be modified such that it satisfies  certain uniform Lipschitzity and boundedness. In Section 3, Theorem \ref{two} and Theorem \ref{two2} will be verified for the modified Lagrangian. According to a priori compactness estimate, the relaxation from the additional assumptions  to the Osgood growth assumption will achieved by a limiting passage in Section 4. Moreover, the proofs of Theorem \ref{two}, Theorem \ref{two2} and Theorem \ref{two1} for the original Lagrangian will be completed.

  \section{\sc Modification of the Lagrangian function}
In this section, in order to avoid the singularity of $h_{x_0,u_0}(x,t)$ at $t=0$ for $x_0\neq x$, we modify the  the Lagrangian function $L(x,u,\dot{x})$ by $L_R(x,u,\dot{x})$ such that $L_R(x,u,\dot{x})$ satisfies:
\begin{itemize}
\item [\textbf{(A1)}]  \textbf{Uniform Lipschitzity}: $L_R(x,u,\dot{x})$ is uniformly Lipschitz with respect to $u$. \vspace{-1.5ex}
\item [\textbf{(A2)}]  \textbf{Uniform boundedness}: $L_R(x,u,\dot{x})-L_R(x,0,\dot{x})$ is uniformly bounded.
%\item [\textbf{(A3)}] \textbf{Completeness of the Flow}: The flows generated by $L_R(x,u,\dot{x})$  are complete;
\end{itemize}
We denote
  \begin{equation}\label{199}
  V(x,u,\dot{x}):=L(x,u,\dot{x})-L(x,0,\dot{x}).
   \end{equation}For a given $R>0$, we choose a $C^\infty$ function $\rho_R(u)$ such that
 \begin{equation}\label{roooo}
\rho_R(u)=\left\{\begin{array}{ll}
\hspace{-0.4em}1,&  |u|\leq R,\\
\hspace{-0.4em}0,&|u|>R+1,\\
\end{array}\right.
\end{equation}
otherwise $0<\rho_R(u)<1$. Without loss of generality, one can require $|\rho'_R(u)|<2$.
%Moreover,
%we denote
%\begin{equation}\label{200}
%V_R(x,u,\dot{x}):=\rho_R(u)V(x,u,\dot{x}).
%\end{equation}It is easy to see that $V_R(x,u,\dot{x})$ is uniformly Lipschitz with respect to $u$. We  denote the Lipschitz constant of $V_R(x,u,\dot{x})$ by $\lambda_R$.
%From (\ref{199}) and (\ref{200}), it follows that  $V_R(x,0,\dot{x})=0$. Hence, $|V_R(x,u,\dot{x})|=0$ for $|u|>R+1$, otherwise, we have
%\begin{equation}\label{br}
%|V_R(x,u,\dot{x})|=|V_R(x,u,\dot{x})-V_R(x,0,\dot{x})|\leq \lambda_R|u|\leq \lambda_R (R+1).
%\end{equation}
 Let
\begin{equation}\label{lr}
\bar{L}_R(x,u,\dot{x})=L(x,0,\dot{x})+\rho_R(u)V(x,u,\dot{x}).
\end{equation}
 We construct a $C^r$ ($r\geq 2$) function denoted by $L_{R}(x,u,\dot{x})$ satisfying
\begin{equation}\label{llr}
L_{R}(x,u,\dot{x}):=\alpha_R (\dot{x})\bar{L}_R(x,u,\dot{x})+\mu_R\beta(|\dot{x}|^2-R^2),
\end{equation}where $\bar{L}_R$ is defined as (\ref{lr}) and $\alpha_R (\dot{x})$ is a $C^\infty$ function satisfying
\begin{equation}\label{alpha}
\alpha_R (\dot{x})=\left\{\begin{array}{ll}
\hspace{-0.4em}1,&  |\dot{x}|\leq R+1,\\
\hspace{-0.4em}0,&|\dot{x}|>R+2,\\
\end{array}\right.
\end{equation}
otherwise $0<\alpha_R (\dot{x})<1$. Without loss of generality, one can require $|\alpha'_R (\dot{x})|<2$ and $|\alpha''_R (\dot{x})|<2$. $\mu_R$ is a sufficient large constant, which is determined by (\ref{murr}) below.
$\beta(z) $ is a $C^\infty$ function satisfying  \begin{equation}\label{beta}
\left\{\begin{array}{ll}
\hspace{-0.4em}\beta(z)=0,&  z\leq 0,\\
\hspace{-0.4em}\beta(z)>z^2-1,&z> 0,\\
\hspace{-0.4em}\beta'(z)>0,&z> 0,\\
\hspace{-0.4em}\beta''(z)>0,&z> 0,\\
\hspace{-0.4em}\beta''(z)>1,&z\geq 1.\\
\end{array}\right.
\end{equation}The existence of $\beta(z)$ can be verified easily. Clearly,  $L_{R}(x,u,\dot{x})$ converges uniformly on compact sets to $L$ as $R\rightarrow\infty$.

It is easy to see that for a given $R>0$, there exists $\lambda_{R}>0$ such that
\begin{equation}
|L_{R}(x,u,\dot{x})-L_{R}(x,v,\dot{x})|\leq\lambda_{R}|u-v|.
\end{equation}
\begin{equation}
|L_{R}(x,u,\dot{x})-L_{R}(x,0,\dot{x})|\leq\lambda_{R}|u|\leq \lambda_{R}(R+1).
\end{equation}
Therefore, $L_{R}(x,u,\dot{x})$ satisfies (A1), (A2).

In addition, $L_{R}(x,u,\dot{x})$  satisfies (L3). In fact,
it suffices to consider the case with $|u|>R+1$ and $|\dot{x}|>R+1$. In this case,
\[L_R(x,u,\dot{x})=\mu_R\beta(|\dot{x}|^2-R^2),\]
which is an integrable system. Thus, $L_{R}(x,u,\dot{x})$  satisfies (L3).

In the following, we show that $L_{R}(x,u,\dot{x})$  satisfies (L1) and (L2).

\vspace{1em}

\noindent\textbf{Claim A:} $L_{R}(x,u,\dot{x})$  satisfies (L1).

\noindent\textbf{Proof of the Claim A:}
It suffices to show that for given $x\in M$ and $u\in \R$,  $\partial^2 L_R/\partial {\dot{x}}^2 (x,u,\dot{x})>0$. We only consider the case with $R<|u|\leq R+1$, the other cases are similar.
\begin{itemize}
\item [(i)] For $|\dot{x}|\leq R$,
\[L_R(x,u,\dot{x})=\bar{L}_R(x,u,\dot{x})=L(x,0,\dot{x})+\rho(u)(L(x,u,\dot{x})-L(x,0,\dot{x})).\]
Hence, we have
\[\frac{\partial^2 L_R}{\partial {\dot{x}}^2} (x,u,\dot{x})=\frac{\partial^2 \bar{L}_R}{\partial {\dot{x}}^2} (x,u,\dot{x})=\rho(u)\frac{\partial^2 L}{\partial {\dot{x}}^2} (x,u,\dot{x})+(1-\rho(u))\frac{\partial^2 L}{\partial {\dot{x}}^2} (x,0,\dot{x})>0.\]
\item [(ii)] For $R<|\dot{x}|\leq R+1$,
\begin{align*}
L_R(x,u,\dot{x})=\bar{L}_R(x,u,\dot{x})+\mu_R\beta(|\dot{x}|^2-R^2).
\end{align*}
It follows that
\[\frac{\partial^2 L_R}{\partial {\dot{x}}^2} (x,u,\dot{x})=\frac{\partial^2 \bar{L}_R}{\partial {\dot{x}}^2} (x,u,\dot{x})+2\mu_R\left(2\beta''(|\dot{x}|^2-R^2)|\dot{x}|^2+\beta'(|\dot{x}|^2-R^2)\right)>0.\]
\item [(iii)] For $R+1<|\dot{x}|\leq R+2$,
\begin{align*}
L_R(x,u,\dot{x})=\alpha_R (\dot{x}) \bar{L}_R(x,u,\dot{x})+\mu_R\beta(|\dot{x}|^2-R^2).
\end{align*}
It yields that
\begin{align*}
\frac{\partial^2 L_R}{\partial {\dot{x}}^2} (x,u,\dot{x})=&\alpha''_R (\dot{x})\bar{L}_R(x,u,\dot{x})+2\alpha'_R (\dot{x})\frac{\partial \bar{L}_R}{\partial {\dot{x}}} (x,u,\dot{x})+
\alpha_R (\dot{x})\frac{\partial^2 \bar{L}_R}{\partial {\dot{x}}^2} (x,u,\dot{x})\\
&+2\mu_R\left(2\beta''(|\dot{x}|^2-R^2)|\dot{x}|^2+\beta'(|\dot{x}|^2-R^2)\right).
\end{align*}
Based on the construction of $\alpha_R (\dot{x})$ and the compactness of $M$, we take
\begin{equation}\label{murr}
\mu_R>\max\left\{\frac{\max|\bar{L}_R(x,u,\dot{x})|+\max |\frac{\partial \bar{L}_R}{\partial {\dot{x}}} (x,u,\dot{x})|
 }{(R+1)^2}, 1\right\}
 \end{equation}
then $\partial^2 L_R/\partial {\dot{x}}^2 (x,u,\dot{x})>0$.
\item [(iv)] For $|\dot{x}|> R+2$,
\[L_R(x,u,\dot{x})=\mu_R\beta(|\dot{x}|^2-R^2),\]
which implies
\[\frac{\partial^2 L_R}{\partial {\dot{x}}^2} (x,u,\dot{x})=2\mu_R\left(2\beta''(|\dot{x}|^2-R^2)|\dot{x}|^2+\beta'(|\dot{x}|^2-R^2)\right)>0\]
\end{itemize}
Therefore, $L_{R}(x,u,\dot{x})$  satisfies (L1).
\End

\vspace{1em}

\noindent\textbf{Claim B:} $L_{R}(x,u,\dot{x})$  satisfies (L2).

\noindent\textbf{Proof of the Claim B:}
Without loss of generality, let $I:=[-C,C]$. By (L2), for $u\in I$, there exists $L_I(x,\dot{x})$ has the property that for any $A>0$, one can find $C_A>0$ such that
\begin{equation}\label{supergr1}
L_I(x,\dot{x})\geq A|\dot{x}|-C_A,
\end{equation}
and
\begin{equation}
L(x,u,\dot{x})\geq L_I(x,\dot{x}).
\end{equation}
In addition, $L_I(x,\dot{x})$ has a lower bound denoted by $D<0$. It follows that
\begin{equation}\label{keyfor1}
\begin{split}
L_R(x,u,\dot{x})&=\alpha_R (\dot{x})\bar{L}_R(x,u,\dot{x})+\mu_R\beta(|\dot{x}|^2-R^2),\\
&=\alpha_R (\dot{x})(\rho(u)L(x,u,\dot{x})+(1-\rho(u))L(x,0,\dot{x}))+\mu_R\beta(|\dot{x}|^2-R^2),\\
&\geq\alpha_R (\dot{x})(\rho(u)L_I(x,\dot{x})+(1-\rho(u))L_I(x,\dot{x}))+\mu_R\beta(|\dot{x}|^2-R^2),\\
&\geq\alpha_R (\dot{x})L_I(x,\dot{x})+\mu_R\beta(|\dot{x}|^2-R^2).
\end{split}
\end{equation}
 According to (\ref{alpha}) and (\ref{beta}), it follows from (\ref{supergr1}) that
\begin{itemize}
\item [(i)] for $|\dot{x}|\leq R$,
\[L_R(x,u,\dot{x})\geq L_I(x,\dot{x})\geq A|\dot{x}|-C_A;\]
\item [(ii)] for $R<|\dot{x}|\leq R+1$,
\begin{align*}
L_R(x,u,\dot{x})&\geq L_I(x,\dot{x})+\mu_R\beta(|\dot{x}|^2-R^2),\\
&\geq L_I(x,\dot{x})\geq A|\dot{x}|-C_A;
\end{align*}
\item [(iii)] for $R+1<|\dot{x}|\leq R+2$,
\begin{align*}
L_R(x,u,\dot{x})&\geq \alpha_R (\dot{x}) L_I(x,\dot{x})+\mu_R\beta(|\dot{x}|^2-R^2),\\
&\geq D+|\dot{x}|^2-1\geq A|\dot{x}|-C'_A;
\end{align*}
\item [(iv)] for $|\dot{x}|> R+2$,
\[L_R(x,u,\dot{x})\geq \mu_R\beta(|\dot{x}|^2-R^2)\geq |\dot{x}|^2-1+D\geq A|\dot{x}|-C'_A.\]
\end{itemize}
Take $C''_A:=\max\{C_A,C'_A\}$. Given $R>0$, for any $A>0$, there exists $C''_A>0$ independent of $R$ such that
\begin{equation}\label{import}
L_R(x,u,\dot{x})\geq A|\dot{x}|-C''_A.
\end{equation}
Therefore, $L_{R}(x,u,\dot{x})$  satisfies (L2).
\End

%\vspace{1em}
%
%\noindent\textbf{Claim C:} $L_{R}(x,u,\dot{x})$  satisfies (L3).
%
%\noindent\textbf{Proof of the Claim C:}
%It suffices to consider the case with $|u|>R+1$ and $|\dot{x}|>R+1$. In this case,
%\[L_R(x,u,\dot{x})=\mu_R\beta(|\dot{x}|^2-R^2),\]
%which is an integrable system. Thus, $L_{R}(x,u,\dot{x})$  satisfies (L3).
%\End

Therefore, $L_{R}(x,u,\dot{x})$ satisfies (L1)-(L3), (A1) and (A2).

For the simplicity of notations, we use $C$ to denote a positive constant only depending on the length of time interval $T$. Without ambiguity, we do not distinguish the constants in different formulas.

\section{\sc Deduction for the modified Lagrangian}
In this section, we are devoted to proving the theorems for the modified Lagrangian (\ref{llr}). Without ambiguity, we still use $L$ to denote $L_{R}$ for a given $R>0$.

The proof of Theorem \ref{two}  will be proceeded  by three steps. In the first step, we will prove the existence and uniqueness  of $h_{x_0,u_0}(x,t)$. Moreover, we will verify a useful property of $h_{x_0,u_0}(x,t)$ so called Markov property. In the second step, we will show that the relation between calibrated curves and solutions of (\ref{hjech}).
%Meanwhile, Theorem \ref{two1} will be verified in this step.
 In the third step, we will  complete the proof of Theorem \ref{two2}.
\subsection{Existence and uniqueness of the action function}
In this step, we are concerned with the existence and uniqueness  of $h_{x_0,u_0}(x,t)$. In order to use the Tonelli theorem to verify the existence of the curve achieving the infimum of (\ref{kkkk}), we consider Theorem \ref{two} in the set of the absolutely continuous curves. To fix the notions,
We use $C^{\text{ac}}([0,t],M)$ to denote the set of all absolutely continuous curves $\gm:[0,t]\rightarrow M$. First of all, we verify the existence of $h_{x_0,u_0}(x,t)$. Unfortunately, classical fixed points theorems can not be used due to the singularity of $h_{x_0,u_0}(x,t)$ at $t=0$ and $x\neq x_0$. Then we have to find a new way to prove that.
\begin{Lemma}\label{exist}
For given $x_0\in M$, $u_0\in\R$, there exists $h_{x_0,u_0}(x,t)\in C(M\times(0,T],\R)$ such that
\begin{equation}\label{u321}
h_{x_0,u_0}(x,t)=u_0+\inf_{\substack{\gm(t)=x \\  \gm(0)=x_0 } }\int_0^tL(\gm(\tau),h_{x_0,u_0}(\gm(\tau),\tau),\dot{\gm}(\tau))d\tau,
\end{equation}where $\gm\in C^{\text{ac}}([0,t],M)$.
\end{Lemma}
\Proof   Without ambiguity, we omit the subscripts $x_0$ and $u_0$ of $h_{x_0,u_0}(x,t)$. We consider a sequence generated by the following iteration formally:
\begin{equation}\label{ui}
h_{i+1}(x,t)=u_0+\inf_{ \gm(t)=x \atop  \gm(0)=x_0 }\int_0^tL(\gm(\tau),h_{i}(\gm(\tau),\tau),\dot{\gm}(\tau))d\tau,
\end{equation}where $i=0,1,2,\ldots$ and $h_0(x,t)=0$.

By the construction of $h_i$, there holds,
\begin{equation}
h_{1}(x,t)=u_0+\inf_{ \gm(t)=x \atop  \gm(0)=x_0 }\int_0^tL(\gm(\tau),0,\dot{\gm}(\tau))d\tau,
\end{equation}
It is easy to see that $h_1(x,t)$ is continuous with respect to $(x,t)\in M\times(0,T]$.
%In order to iterate, we set
%\[h_{i}(x,0)=u_0,\quad i=1,2,\ldots.\]
Let
\[V(x,u,\dot{x}):=L(x,u,\dot{x})-L(x,0,\dot{x}).\]
Then we have
\begin{equation}\label{uiii}
\begin{split}
h_{2}(x,t)&=u_0+\inf_{ \gm(t)=x \atop  \gm(0)=x_0 }\int_0^tL(\gm(\tau),h_{1}(\gm(\tau),\tau),\dot{\gm}(\tau))d\tau,\\
&=u_0+\inf_{ \gm(t)=x \atop  \gm(0)=x_0}\int_0^tL(\gm(\tau),0,\dot{\gm}(\tau))+V(\gm(\tau),h_{1}(\gm(\tau),\tau),\dot{\gm}(\tau))d\tau.
\end{split}
\end{equation}

\noindent\textbf{Claim:} There exists an absolutely continuous curve $\gm_2: [0,t]\rightarrow M$ satisfying $\gm_2(0)=x_0$  and $\gm_2(t)=x$ such that the infimum in (\ref{uiii}) can be achieved.

\noindent\textbf{Proof of the Claim:}
Let
\[A_{h_1}(\gm|_{[0,t]}):=\int_0^tL(\gm(\tau),h_{1}(\gm(\tau),\tau),\dot{\gm}(\tau))d\tau,\quad A_{0}(\gm|_{[0,t]}):=\int_0^tL(\gm(\tau),0,\dot{\gm}(\tau))d\tau.\]

By (A2), we have
\begin{equation}\label{lipkkey}
|A_{h_1}(\gm|_{[0,t]})-A_{0}(\gm|_{[0,t]})|\leq \int_0^t |V(\gm(\tau),h_{1}(\gm(\tau),\tau),\dot{\gm}(\tau))|d\tau\leq Ct,
\end{equation}
where $C:=\sup_{(x,u,\dot{x})\in TM\times\R}|V(x,u,\dot{x})|$.  Based on the compactness of $M$,
it follows from Tonelli's theorem  \cite{M1} that for $K\in \R$, the set
\[\Sigma_K:=\{\gm\in C^{\text{ac}}([0,t],M)\  | \ A_{0}(\gm|_{[0,t]})\leq K\}\]
 is a compact subset of $C^{\text{ac}}([0,t],M)$ endowed with the topology of uniform convergence. Based on (\ref{lipkkey}), it is easy to see that
 for $K'\in \R$, the set
\[\Sigma_{K'}:=\{\gm\in C^{\text{ac}}([0,t],M)\  | \ A_{h_1}(\gm|_{[0,t]})\leq K'\}\]
 is a precompact subset of $C^{\text{ac}}([0,t],M)$ endowed with the topology of uniform convergence.

 Let
 \[L(x,\dot{x},\tau):=L(x,h_1(x,\tau),\dot{x})=L(x,0,\dot{x})+V(x,h_1(x,\tau),\dot{x}).\]
It is easy to see that  $L(x,\dot{x},\tau)$  and $\partial L(x,\dot{x},\tau)/\partial\dot{x}$ are continuous. Hence, a standard argument yields that $A_{h_1}$ is lower semi-continuous, which implies $\Sigma_{K'}$ is compact (see  \cite{F3} Section 3.2 for instance).

Let \[C_{\text{inf}}=\inf_{ \gm(t)=x \atop  \gm(0)=x_0 } A_{h_1}(\gm|_{[0,t]}),\]
where $\gm\in C^{\text{ac}}([0,t],M)$. For each integer $n\geq 1$, we denote
\[\mathcal{C}_n:=\left\{\gm\in C^{\text{ac}}([0,t],M)\ \bigg|\ \gm(0)=x_0, \gm(t)=x, A_{h_1}(\gm|_{[0,t]})\leq C_{\text{inf}}+\frac{1}{n}.\right\}\]
Since $\Sigma_{K'}$ is compact for any $K'\in \R$, the set is a nonempty compact subset of $C^{\text{ac}}([0,t],M)$.
$\mathcal{C}_n$ is decreasing, then $\cap_{n\geq 1}\mathcal{C}_n$ is nonempty. Taking $\gm_2\in \cap_{n\geq 1}\mathcal{C}_n$,  we have $\gm_2(0)=x_0$, $\gm_2(t)=x$ and
\[A_{h_1}(\gm_2|_{[0,t]})=C_{\text{inf}}.\]
 Hence, the infimum in (\ref{uiii}) can be achieved at an absolutely continuous curve $\gm_2: [0,t]\rightarrow M$ satisfying $\gm_2(0)=x_0$  and $\gm_2(t)=x$. \End

 Moreover, $h_{2}(x,t)$ is continuous with respect to $(x,t)\in M\times(0,T]$.
Inductively, for each $i>2$, $h_{i}(x,t)$ is continuous with respect to $(x,t)\in M\times(0,T]$ and there exists an absolutely continuous curve $\gm_i: [0,t]\rightarrow M$ satisfying $\gm_i(0)=x_0$  and $\gm_i(t)=x$ such that the infimum in (\ref{ui}) can be achieved. To fix the notions, $\gm_i$ is called a minimal curve of $h_i$.

For any $x\in M$, let $\gm_2:[0,t]\rightarrow M$ be a minimal curve of $h_2$ with $\gm_2(0)=x_0$ and  $\gm_2(t)=x$, $\gm_1:[0,s]\rightarrow M$ be a minimal curve of $h_1$ with $\gm_1(0)=x_0$ and  $\gm_1(s)=\gm_2(s)$. By (A1), we have
\begin{equation}\label{lip}
|L(x,u,\dot{x})-L(x,v,\dot{x})|\leq \lambda|u-v|.
\end{equation}Then for $s\in [0,t]$, we have
\begin{align*}
&\ \ \ h_2(\gm_2(s),s)-h_1(\gm_2(s),s)\\
&\leq\int_0^s L(\gm_1(\tau),h_1(\gm_1(\tau),\tau),\dot{\gm}_1(\tau))d\tau-\int_0^sL(\gm_1(\tau),0,\dot{\gm}_1(\tau))d\tau,\\
&\leq\int_0^s |L(\gm_1(\tau),h_1(\gm_1(\tau),\tau),\dot{\gm}_1(\tau))-L(\gm_1(\tau),0,\dot{\gm}_1(\tau))|d\tau,
\end{align*}
which together with (A2) implies
\begin{equation}\label{2kt}
h_2(\gm_2(s),s)-h_1(\gm_2(s),s)\leq Cs,
\end{equation} By a similar argument, we have
\begin{equation}
 h_2(\gm_2(s),s)-h_1(\gm_2(s),s)\geq -Cs.
\end{equation}
Hence, there holds that for $s\in [0,t]$,
\begin{equation}
 |h_2(\gm_2(s),s)-h_1(\gm_2(s),s)|\leq Cs.
\end{equation}
In particular, there holds for  $(x,t)\in M\times (0,T]$,
\[\|h_2(x,t)-h_1(x,t)\|_{C^0}\leq CT,\]
where $\|\cdot\|_{C^0}$ denotes the $C^0$-norm of continuous function space defined on $M\times (0,T]$.

 Let $\gm_3:[0,t]\rightarrow M$ be a minimal curve of $h_3$ with $\gm_3(0)=x_0$ and  $\gm_3(t)=x$,  $\gm_2:[0,s]\rightarrow M$ be a minimal curve of $h_2$ with $\gm_2(0)=x_0$ and  $\gm_2(s)=\gm_3(s)$. Moreover,
\begin{align*}
&\ \ \ h_3(\gm_3(s),s)-h_2(\gm_3(s),s)\\
&\leq\int_0^s L(\gm_2(\tau),h_2(\gm_2(\tau),\tau),\dot{\gm}_2(\tau))d\tau-\int_0^s L(\gm_2(\tau),h_1(\gm_2(\tau),\tau),\dot{\gm}_2(\tau))d\tau,\\
&\leq\int_0^s |L(\gm_2(\tau),h_2(\gm_2(\tau),\tau),\dot{\gm}_2(\tau))-L(\gm_2(\tau),h_1(\gm_2(\tau),\tau),\dot{\gm}_2(\tau))|d\tau,\\
&\leq\lambda\int_0^s |h_2(\gm_2(\tau),\tau))-h_1(\gm_2(\tau),\tau)|d\tau\leq \lambda C\int_0^s \tau d\tau=\frac{1}{2}C\lambda s^2.
\end{align*}
 By a similar argument, we have
 \begin{equation}
h_3(\gm_3(s),s)-h_2(\gm_3(s),s)\geq-\frac{1}{2}C\lambda s^2.
\end{equation}
In particular, we have
\begin{equation}\label{C32}
\|h_3(x,t)-h_2(x,t)\|_{C^0}\leq\frac{1}{2}C\lambda T^2,
\end{equation}
Repeating the argument above $n$ times, we have
\begin{equation}\label{Cn}
\|h_{n+1}(x,t)-h_n(x,t)\|_{C^0}\leq \frac{1}{n!}C\lambda^{n-1} T^n.
\end{equation}It follows from (\ref{Cn}) that as $n\rightarrow \infty$,
\begin{equation}
\|h_{n+1}(x,t)-h_n(x,t)\|_{C^0}\rightarrow 0,
\end{equation}
hence there exists $\bar{h}(x,t)\in C(M\times (0,T],\R)$ such that
\begin{equation}
\lim_{n\rightarrow\infty}h_n(x,t)=h_1(x,t)+\lim_{n\rightarrow\infty}\sum_{i=1}^{n-1}\left(h_{i+1}(x,t)-h_i(x,t)\right)=\bar{h}(x,t),
\end{equation}where the limiting passage is uniform on $M\times (0,T]$ and $\bar{h}(x,t)$ satisfies (\ref{u321}). This finishes the proof of Lemma \ref{exist}.
\End

Lemma \ref{exist} implies that there exists $h_{x_0,u_0}(x,t)\in C(M\times(0,T],\R)$ such that
\begin{equation}\label{hh}
h_{x_0,u_0}(x,t)=u_0+\inf_{\substack{\gm(t)=x \\  \gm(0)=x_0 } }\int_0^tL(\gm(\tau),h_{x_0,u_0}(\gm(\tau),\tau),\dot{\gm}(\tau))d\tau.
\end{equation} In particular, the infimum can be achieved at an absolutely continuous curve denoted by $\bar{\gm}$. In the way analogous to Mather theory \cite{M1,Mvc},  we call $h_{x_0,u_0}(x,t)$ an {\it action function} with respect to $L(x,u,\dot{x})$. To fix the notions, we call $\bar{\gm}$  a {\it calibrated curve} (at finite time).

The following lemma implies the uniqueness of $h_{x_0,u_0}(x,t)$.

\begin{Lemma}\label{unique}
If both $h_{x_0,u_0}(x,t)$ and $g_{x_0,u_0}(x,t)$ satisfy (\ref{hh}), then for $(x,t)\in M\times (0,T]$, we have \[h_{x_0,u_0}(x,t)=g_{x_0,u_0}(x,t).\]
\end{Lemma}
The proof of Lemma \ref{unique} depends on an inequality as follows.

\noindent\textbf{Gronwall's inequality:} Let $F:[0,t]\rightarrow\R$ be continuous and nonnegative. Suppose $C\geq 0$ and $\lambda\geq 0$ are such that for any $s\in [0,t]$,
 \begin{equation}
F(s)\leq C+\int_{0}^s\lambda F(\tau)d\tau.
\end{equation}Then, for any $s\in [0,t]$,
 \begin{equation}
F(s)\leq Ce^{\lambda s}.
\end{equation}
Taking $C=0$, we have the following lemma.
\begin{Lemma}\label{gron}
Let $F:[0,t]\rightarrow\R$ be continuous with $F(0)=0$ and $F(s)>0$ for $s\in (0,t]$. Then for a given $\lambda\geq 0$, there exists $s_0\in (0,t]$ such that
 \begin{equation}\label{Ft}
F(s_0)>\int_{0}^{s_0}\lambda F(\tau)d\tau.
\end{equation}
\end{Lemma}

\noindent\textbf{Proof of Lemma \ref{unique}:} The same as the notations in Lemma \ref{exist}, we denote $h_{x_0,u_0}(x,t)$ and $g_{x_0,u_0}(x,t)$ by $h(x,t)$ and $g(x,t)$ respectively.

On one hand, we will prove
\begin{equation}
h(x,t)\leq g(x,t).
\end{equation}
By contradiction, we assume $h(x,t)> g(x,t)$.
 Let $\gm_g$ be a calibrated curve of $g$ with $\gm_g(0)=x_0$, $\gm_g(t)=x$. We denote
 \begin{equation}
F(\tau):=h(\gm_g(\tau),\tau)-g(\gm_g(\tau),\tau),
\end{equation}where $\tau\in [0,t]$. By (\ref{hh}), we have $F(0)=0$. The assumption $h(x,t)> g(x,t)$ implies $F(t)>0$. Hence, there exists $\tau_0\in [0,t)$ such that $F(\tau_0)=0$ and  $F(\tau)> 0$ for $\tau> \tau_0$.
Let $\gm_h$ be a calibrated curve of $h$ with $\gm_h(0)=x_0$, $\gm_h(\tau_0)=\gm_g(\tau_0)$. For $s\in [\tau_0,t]$, we construct a family of $\gm_s:[0,s]\rightarrow M$ as follows:
\begin{equation}
\gm_s(\tau)=\left\{\begin{array}{ll}
\hspace{-0.4em}\gm_h(\tau),& \tau\in [0,\tau_0],\\
\hspace{-0.4em}\gm_g(\tau),&\tau\in (\tau_0,s].\\
\end{array}\right.
\end{equation}
 Based on the definition of $h(x,t)$ (see (\ref{hh})), we have
 \begin{align*}
 h(\gm_g(s),s)&=u_0+\inf_{\substack{\gm(s)=\gm_g(s) \\  \gm(0)=x_0 } }\int_0^sL(\gm(\tau),h(\gm(\tau),\tau),\dot{\gm}(\tau))d\tau,\\
 &\leq u_0+\int_{0}^sL(\gm_s(\tau),h(\gm_s(\tau),\tau),\dot{\gm}_s(\tau))d\tau,\\
 &= h(\gm_g(\tau_0),\tau_0)+\int_{\tau_0}^sL(\gm_g(\tau),h(\gm_g(\tau),\tau),\dot{\gm}_g(\tau))d\tau.
 \end{align*}
Similarly, for $g(x,t)$, we have
 \begin{equation}
g(\gm_g(s),s)=g(\gm_g(\tau_0),\tau_0)+\int_{\tau_0}^sL(\gm_g(\tau),g(\gm_g(\tau),\tau),\dot{\gm}_g(\tau))d\tau.
\end{equation}
Since $h(\gm_g(\tau_0),\tau_0)-g(\gm_g(\tau_0),\tau_0)=F(\tau_0)=0$,   then we have
 \begin{align*}
 &h(\gm_g(s),s)-g(\gm_g(s),s)\\
 &\leq \int_{\tau_0}^sL(\gm_g(\tau),h(\gm_g(\tau),\tau),\dot{\gm}_g(\tau))-
 L(\gm_g(\tau),g(\gm_g(\tau),\tau),\dot{\gm}_g(\tau))d\tau,\\
 &\leq \int_{\tau_0}^s\lambda |h(\gm_g(\tau),\tau)-g(\gm_g(\tau),\tau)|d\tau,\\
 &=\int_{\tau_0}^s\lambda (h(\gm_g(\tau),\tau)-g(\gm_g(\tau),\tau))d\tau,
 \end{align*}where the second inequality is owing to (A1).
It follows that for any $s\in (\tau_0,t]$
 \begin{equation}\label{Ft1}
F(s)\leq \int_{\tau_0}^s\lambda F(\tau)d\tau,
\end{equation}
which is in contradiction with Lemma \ref{gron}.  Thus, we obtain $h(x,t)\leq g(x,t)$.

On the other hand, it follows from a similar argument that $h(x,t)\geq g(x,t)$. So far, we have shown that
$h(x,t)=g(x,t)$ for $(x,t)\in M\times (0,T]$, which finishes the proof of Lemma \ref{unique}.
\End

 Lemma \ref{exist} and Lemma \ref{unique} imply the well posedness of $h_{x_0,u_0}(x,t)$. There exists a unique $h_{x_0,u_0}(x,t)\in C(M\times(0,T],\R)$ such that
\begin{equation}
h_{x_0,u_0}(x,t)=u_0+\inf_{\substack{\gm(t)=x \\  \gm(0)=x_0 } }\int_0^tL(\gm(\tau),h_{x_0,u_0}(\gm(\tau),\tau),\dot{\gm}(\tau))d\tau,
\end{equation}where $\gm\in C^{\text{ac}}([0,t],M)$. Similarly, it is easy to see that there exists a unique $h^+_{x_0,u_0}(x,t)\in C(M\times[-T,0),\R)$ such that
\begin{equation}
h^+_{x_0,u_0}(x,t)=u_0-\inf_{\substack{ \gm(0)=x_0 \\ \gm(t)=x } }\int_t^0L(\gm(\tau),h^+_{x_0,u_0}(\gm(\tau),\tau),\dot{\gm}(\tau))d\tau,
\end{equation}

We only discuss $h_{x_0,u_0}(x,t)$ in the following. The corresponding properties of $h^+_{x_0,u_0}(x,t)$ can be obtained similarly.

Based on a similar argument as Lemma \ref{unique}, we prove that $h_{x_0,u_0}(x,t)$ satisfies Markov property, which is useful in  further research.
% It insures that the differentiate point  can be chose arbitrarily in the proof of Lemma \ref{chc} below, from which the coincidence between $h_{x_0,u_0}(x,t)$ and $U(t)$ can be extended globally.
 We omit the subscripts $x_0$ and $u_0$ of $h_{x_0,u_0}(x,t)$ for  simplicity.

 \begin{Lemma}\label{tria1}
  \begin{equation}\label{trihx}
h(x,t+s)=\inf_{y\in M}h_{y,h(y,t)}(x,s).
\end{equation}
 \end{Lemma}
\Proof
On one hand, we will prove $h(x,t+s)\geq\inf_{y\in M}h_{y,h(y,t)}(x,s)$. Let $\gm_1:[0,t+s]\rightarrow M$ be a calibrated curve of $h$ with $\gm_1(0)=x_0$ and $\gm_1(t+s)=x$. Consider $\bar{y}\in \gm_1$ with $\gm_1(t)=\bar{y}$. It suffices to show
  \begin{equation}\label{s11}
h(x,t+s)\geq h_{\bar{y},h(\bar{y},t)}(x,s).
\end{equation}
By contradiction, we assume $h(x,t+s)< h_{\bar{y},h(\bar{y},t)}(x,s)$. By the definition of $h(x,t+s)$, we have
\begin{align*}
h(x,t+s)&=u_0+\int_0^{t+s}L(\gm_1(\tau),h(\gm_1(\tau),\tau),\dot{\gm}_1(\tau))d\tau,\\
&=h(\bar{y},t)+\int_t^{t+s}L(\gm_1(\tau),h(\gm_1(\tau),\tau),\dot{\gm}_1(\tau))d\tau,\\
&=h(\bar{y},t)+\int_0^{s}L(\gm_1(\sigma+t),h(\gm_1(\sigma+t),\sigma+t),\dot{\gm}_1(\sigma+t))d\sigma.
\end{align*}
By the definition of $h_{\bar{y},h(\bar{y},t)}(x,s)$, we have
\begin{align*}
h_{\bar{y},h(\bar{y},t)}(x,s)&=h(\bar{y},t)+\inf_{\substack{\gm(s)=x \\  \gm(0)=\bar{y} } }
\int_0^{s}L(\gm(\tau),h_{\bar{y},h(\bar{y},t)}(\gm(\tau),\tau),\dot{\gm}(\tau))d\tau,\\
&\leq h(\bar{y},t)+\int_0^{s}L(\gm_1(\sigma+t),h_{\bar{y},h(\bar{y},t)}(\gm_1(\sigma+t),\sigma),\dot{\gm}(\sigma+t))d\sigma.
\end{align*}
We denote $u(\sigma):=h(\gm_1(\sigma+t),\sigma+t)$ and $v(\sigma):=h_{\bar{y},h(\bar{y},t)}(\gm_1(\sigma+t),\sigma)$. In particular, we have $h(x,t+s)=u(s)$ and $h_{\bar{y},h(\bar{y},t)}(x,s)=v(s)$. Let
 \begin{equation}
F(\sigma):=v(\sigma)-u(\sigma),
\end{equation}where $\sigma\in [0,s]$. It is easy to see that $u(0)=h(\bar{y},t)=v(0)$. Then we have $F(0)=0$. The assumption $h(x,t+s)< h_{\bar{y},h(\bar{y},t)}(x,s)$ implies $F(s)>0$. Hence, there exists $\sigma_0\in [0,s)$ such that $F(\sigma_0)=0$ and  $F(\sigma)> 0$ for $\sigma> \sigma_0$. Moreover, for any $\tau\in (\sigma_0,s]$, we have
  \begin{equation}
u(\tau)=u(\sigma_0)+\int_{\sigma_0}^{\tau}L(\gm_1(\sigma+t),u(\sigma),\dot{\gm}_1(\sigma+t))d\sigma.
\end{equation}
Let $\gm_2$ be a calibrated curve of $h$ with $\gm_2(0)=\bar{y}$, $\gm_2(\sigma_0)=\gm_1(\sigma_0+t)$. For $\sigma\in [\sigma_0,\tau]$, we construct a family of $\gm_\tau:[0,\tau]\rightarrow M$ as follows:
\begin{equation}
\gm_\tau(\sigma)=\left\{\begin{array}{ll}
\hspace{-0.4em}\gm_2(\sigma),& \sigma\in [0,\sigma_0],\\
\hspace{-0.4em}\gm_1(\sigma+t),&\sigma\in (\sigma_0,\tau].\\
\end{array}\right.
\end{equation} Moreover,
for any $\tau\in (0,s]$, we have
  \begin{align*}
v(\tau)&\leq h(\bar{y},t)+\int_0^{\tau}L(\gm_\tau(\sigma+t),
h_{\bar{y},h(\bar{y},t)}(\gm_\tau(\sigma),\sigma),\dot{\gm}_\tau(\sigma))d\sigma,\\
&=v(\sigma_0)+\int_{\sigma_0}^{\tau}L(\gm_1(\sigma+t),v(\sigma),\dot{\gm}_1(\sigma+t))d\sigma.
\end{align*}
Since $v(\sigma_0)-u(\sigma_0)=F(\sigma_0)=0$, it follows from (A1) that
  \begin{equation}
v(\tau)-u(\tau)\leq \int_{\sigma_0}^\tau\lambda(v(\sigma)-u(\sigma))d\sigma.
\end{equation}
It yields that for any $\tau\in (\sigma_0,s]$,
 \begin{equation}
F(\tau)\leq \int_{\sigma_0}^\tau\lambda F(\sigma)d\sigma,
\end{equation}
which is in contradiction with Lemma \ref{gron}. Hence,
 \begin{equation}\label{s22}
h(x,t+s)\geq\inf_{y\in M}h_{y,h(y,t)}(x,s).
\end{equation}

On the other hand, we will prove $h(x,t+s)\leq\inf_{y\in M}h_{y,h(y,t)}(x,s)$. Due to the compactness of $M$, there exists $\tilde{y}$ such that $\inf_{y\in M}h_{y,h(y,t)}(x,s)=h_{\tilde{y},h(\tilde{y},t)}(x,s)$.
 Let $\gm_3:[0,t]\rightarrow M$ be a calibrated curve of $h$ with $\gm_3(0)=x_0$ and $\gm_3(t)=\tilde{y}$.  Let $\gm_4:[0,s]\rightarrow M$ be a calibrated curve of $h_{\tilde{y},h(\tilde{y},t)}$ with $\gm_4(0)=\tilde{y}$ and $\gm_4(s)=x$. By a similar argument as (\ref{s11})-(\ref{s22}), we have
 \begin{equation}\label{s23}
h(x,t+s)\leq\inf_{y\in M}h_{y,h(y,t)}(x,s),
\end{equation}
which together with (\ref{s22}) implies
  \begin{equation}
h(x,t+s)=\inf_{y\in M}h_{y,h(y,t)}(x,s).
\end{equation}
This completes the proof of Lemma \ref{tria1}.
\End

From Lemma \ref{tria1}, we obtain
  \begin{equation}\label{triin}
h(x,t+s)-u_0\leq (h_{y,h(y,t)}(x,s)-h(y,t))+(h(y,t)-u_0).
\end{equation}
Let
\[B^t(x,u;y):=h_{x,u}(y,t)-u.\]
Then (\ref{triin}) can be formulated as
\[B^{t+s}(x_0,u_0;x)\leq B^s\left(y,h_{x_0,u_0}(y,t);x\right)+B^t(x_0,u_0;y),\]
which can be seen as a triangle inequality from the dynamical point of view. In particular, the equality holds if and only if $y$ belongs to the calibrated curve $\gm$ of $h$ with $\gm(0)=x_0$ and $\gm(t+s)=x$.
\begin{Remark}
(\ref{trihx}) has different names in different fields. For instance, it is called Dynamical Programming Principle in optimal control theory \cite{Ba3,CS}. It also appears as Semigroup Property, Markov Property, Shortest Path Algorithm in PDE, stochastic control and graph theory respectively \cite{CHL2,S}. Roughly speaking, (\ref{trihx}) implies the connection between the local optimum and the global optimum.
\end{Remark}
\subsection{Calibrated curve and solution}
In the previous step, we obtain that there exists a unique  $h(x,t)\in C(M\times(0,T],\R)$ satisfying (\ref{hh}). In this step, we will show that the relation between calibrated curves and solutions  of (\ref{hjech}). More precisely, we have the following lemma:

\begin{Lemma}\label{chc}
Let $\bar{\gm}:[0,t]\rightarrow M$ be a calibrated curve of $h$, then $\bar{\gm}$ is $C^1$ and for $\tau\in (0,t)$,
$(\bar{\gm}(\tau),u(\tau),p(\tau))$ satisfies the contact Hamiltonian equation (\ref{hjech}),
 where
\begin{equation}\label{ccure}
u(\tau)=h(\bar{\gm}(\tau),\tau)\quad\text{and}\quad p(\tau)=\frac{\partial L}{\partial \dot{x}}(\bar{\gm}(\tau),h(\bar{\gm}(\tau),\tau),\dot{\bar{\gm}}(\tau)).
\end{equation}
\end{Lemma}

\Proof
Since $\bar{\gm}\in C^{ac}([0,t],M)$, then the derivative $\dot{\bar{\gm}}(\tau)$ exists almost everywhere for $\tau\in [0,t]$. Let $t_0\in (0,t)$ be a differentiate point of $\bar{\gm}(\tau)$. It suffices to consider $0<t_0<t$. Denote
\begin{equation}
(x_0,u_0,v_0):=(\bar{\gm}(t_0),h(\bar{\gm}(t_0),t_0),\dot{\bar{\gm}}(t_0)).
 \end{equation}

The proof is divided into three steps. Note that (\ref{hjech}) is the characteristic equation of the Hamilton-Jacobi equation:
\begin{equation}\label{ssshj}
\partial_\tau S(x,\tau)+H(x,S(x,\tau),\partial_x S(x,\tau))=0.
\end{equation}

\noindent\textbf{(a) Construction of the classical solution}

First of all, we will construct a classical solution of (\ref{ssshj}) on a cone-like region (see (\ref{tri}) below). Let $k:=|v_0|$ and
 \[B(0,2k):=\{v:|v|< 2k,\ v\in T_{x_0}M\}.\]
We use $B^*(0,2k)$ to denote the image of $B(0,2k)$ via the Legendre transformation $\mathcal{L}^{-1}:TM\rightarrow T^*M$. That is
\[B^*(0,2k):=\left\{p: p=\frac{\partial L}{\partial v}(x_0,u_0,v),\ v\in B(0,2k)\right\}.\]
  Let $\Psi_\tau:T^*M\times\R\rightarrow T^*M\times\R$ denote the follow generated by the contact Hamiltonian equation (\ref{hjech}). Let $\pi$ be a projection from $T^*M\times\R$ to $ T^*M$ via $(x,u,p)\rightarrow (x,p)$ and let
\[\mathcal{B}_\tau^*(0,2k):=\pi\circ\Psi_{\tau-t_0}(x_0,u_0,B^*(0,2k)).\]Moreover, we denote
\[\Pi_\tau:\mathcal{B}_\tau^*(0,2k)\rightarrow M,\]
%where $\Pi_\tau$ is a projection from $T^*M$ to $M$.
Since the Legendre transformation $\mathcal{L}$ is a diffeomorphism from $T_{x_0}^*M\times\{u_0\}$ to $T_{x}^*M\times\{u_0\}$, it follows that there exists $\epsilon>0$ such that for $\tau\in (t_0,t_0+\epsilon]$,  $\Pi_\tau$ is a diffeomorphism onto its image.

We denote
\[\Omega_\tau:=\Pi_\tau( \mathcal{B}_\tau^*(0,2k)).\]
 We use $\Omega_\epsilon$ to denote the following cone-like region:
\begin{equation}\label{tri}
\Omega_\epsilon:=\{(\tau,x): \tau\in (t_0,t_0+\epsilon),\ x\in \Omega_\tau\}.
 \end{equation}Then for any $(\tau,x)\in\Omega_\epsilon$, there exists a unique $p_0\in B^*(0,2k)$ such that $X(\tau)=x$ where
 \[(X(\tau), U(\tau),P(\tau)):=\Phi_\tau(x_0,u_0,p_0).\]
 Hence, for  any $(\tau,x)\in \Omega_\epsilon$, one can define a $C^1$ function by $S(x,\tau)=U(\tau)$. In particular, we have $S(x,t_0)=u_0$. Moreover, it follows from  the method of characteristics (see \cite{Ar,F3} for instance) that $S(x,\tau)$
is a classical solution of the Hamilton-Jacobi equation (\ref{ssshj}) with $S(x_0,t_0)=u_0$.

\vspace{1em}

\noindent\textbf{(b) Local coincidence between $\bar{\gm}(\tau)$ and $X(\tau)$}

Fix $\tau\in (t_0,t_0+\epsilon)$ and let $S_\tau(x):=S(x,\tau)$. We denote
\begin{equation}
\text{grad}_LS_\tau(x):=\frac{\partial H}{\partial p}(x,S_\tau(x),p),
\end{equation}where $p=\partial_x S_\tau(x)$. In particular, we have $v_0=\text{grad}_LS_{t_0}(x_0)$. It is easy to see that $\text{grad}_LS_\tau(x)$ gives rise to a vector field on $M$.
%Let $\Omega_\epsilon$ be the Legendre transformation of $\Omega_\epsilon$. It is easy to see that
%\[\Omega_\epsilon=\Omega_\epsilon.\]

\noindent\textbf{Claim A:} Let $\gm$ be an absolutely continuous curve with $(\tau,\gm(\tau))\in \Omega_\epsilon$ for  $\tau\in [a,b]\subset [t_0,t_0+\epsilon]$, we have
\begin{equation}\label{hjee}
S(\gm(b),b)-S(\gm(a),a)\leq\int_a^bL(\gm(\tau),
S(\gm(\tau),\tau),\dot{\gm}(\tau))d\tau,
\end{equation}
where the equality holds if and only if $\gm$ is a trajectory of the vector field $\text{grad}_LS_\tau(x)$.

\noindent\textbf{Proof of the claim A:}
From the regularity of $S(x,\tau)$, it follows that
\begin{equation}\label{sasb}
S(\gm(b),b)-S(\gm(a),a)=\int_a^b\left\{\frac{\partial S}{\partial t}(\gm(\tau),\tau)+\langle \frac{\partial S}{\partial x}(\gm(\tau),\tau),\dot{\gm}(\tau)\rangle\right\}d\tau.
\end{equation}By virtue of Fenchel inequality, for each $\tau$ where $\dot{\gm}(\tau)$ exists, we have
\begin{align*}
\left\langle \frac{\partial S}{\partial x}(\gm(\tau),\tau),\dot{\gm}(\tau)\right\rangle\leq &H(\gm(\tau),S(\gm(\tau),\tau),\frac{\partial S}{\partial x}(\gm(\tau),\tau))\\
&+L(\gm(\tau),
S(\gm(\tau),\tau),\dot{\gm}(\tau)).
\end{align*}It follows from (\ref{ssshj}) that for almost every $\tau\in [a,b]$
\begin{equation}
\frac{\partial S}{\partial t}(\gm(\tau),\tau)+\left\langle \frac{\partial S}{\partial x}(\gm(\tau),\tau),\dot{\gm}(\tau)\right\rangle\leq L(\gm(\tau),
S(\gm(\tau),\tau),\dot{\gm}(\tau)).
\end{equation}
By integration, it follows from (\ref{sasb}) that
\begin{equation}\label{ss}
S(\gm(b),b)-S(\gm(a),a)\leq\int_a^bL(\gm(\tau),
S(\gm(\tau),\tau),\dot{\gm}(\tau))d\tau.
\end{equation}We have equality in (\ref{ss}) if and only if the equality holds in the Fenchel inequality, i.e. $\dot{\gm}(\tau)=\text{grad}_LS_\tau(x)$ which means that $\gm$ is a trajectory of the vector field $\text{grad}_LS_\tau(x)$.
\End

\noindent\textbf{Claim B:} For any $\tau\in [t_0,t_0+\epsilon]$, there holds
\begin{equation}
S(\bar{\gm}(\tau),\tau)=h(\bar{\gm}(\tau),\tau),
\end{equation}where $\bar{\gm}$ is a calibrated curve of $h$ with $\bar{\gm}(t_0)=x_0$.

\noindent\textbf{Proof of the claim B:}
This claim can be verified by a similar argument as the one in the proof of Lemma \ref{unique}. A slight difference is that $S(x,t)$ is just well defined for $(t,x)\in \Omega_\epsilon$.

Based on the construction of $\Omega_\epsilon$,  there exists a trajectory of the vector field $\text{grad}_LS_\tau(x)$ denoted by $\tilde{\gm}$ with $\tilde{\gm}(t_0)=x_0$,  $\tilde{\gm}(s)=\bar{\gm}(s)$ and $(\tau, \tilde{\gm}(\tau)\in \Omega_\epsilon$ for any $s\in (t_0,t_0+\epsilon]$. By Claim A, we have
\begin{equation}
\int_{t_0}^{s}L(\tilde{\gm}(\tau),
S(\tilde{\gm}(\tau),\tau),\dot{\tilde{\gm}}(\tau))d\tau=S(\bar{\gm}(s),s)-S(x_0,t_0)\leq\int_{t_0}^{s}L(\gm(\tau),
S(\gm(\tau),\tau),\dot{\gm}(\tau))d\tau,
\end{equation}
where $\gm$ denotes an absolutely continuous curve with $(\tau,\gm(\tau))\in \Omega_\epsilon$.

It is easy to see that the argument  in  the proof of Lemma \ref{unique} only depends on $\bar{\gm}$ and $\tilde{\gm}$. Since both of $\bar{\gm}$ and $\tilde{\gm}$ are contained in $\Omega_\epsilon$, then this claim can be verified by a similar argument as the one in the proof of Lemma \ref{unique}.
 \End

  From the definition of $h_{x_0,u_0}$ (see (\ref{hh})), it follows that
\begin{align*}
S(\bar{\gm}(t_0+\epsilon),t_0+\epsilon)=S(\bar{\gm}(t_0),t_0)
+\int_{t_0}^{t_0+\epsilon}L(\bar{\gm}(\tau),S(\bar{\gm}(\tau),\tau),\dot{\bar{\gm}}(\tau))d\tau,
\end{align*}which together with the claim implies $\bar{\gm}(\tau)$  is a trajectory of the vector field $\text{grad}_LS_\tau(x)$. Let
\begin{equation}
u(\tau):=h(\bar{\gm}(\tau),\tau)\quad\text{and}\quad p(\tau):=\frac{\partial L}{\partial \dot{x}}(\bar{\gm}(\tau),h(\bar{\gm}(\tau),\tau),\dot{\bar{\gm}}(\tau)).
\end{equation}It follows that for $\tau\in [t_0,t_0+\epsilon]$, $(\bar{\gm}(\tau),u(\tau),p(\tau))$  is $C^1$ and satisfies the contact Hamiltonian equation (\ref{hjech}).
%That is
%\[\bar{\gm}(\tau)=X(\tau),\quad u(\tau)=U(\tau),\quad p(\tau)=P(\tau).\]

\vspace{1em}

\noindent\textbf{(c) Global coincidence between $\bar{\gm}(\tau)$ and $X(\tau)$}

We will show that $(\bar{\gm}(\tau),u(\tau),p(\tau))$ coincides with $(X(\tau),U(\tau),P(\tau))$ at $(0,t)$. By the density of the differentiate points of $\bar{\gm}$, it suffices to prove the coincidence at $[t_0,t)$.

 By contradiction, we assume that $(\bar{\gm}(\tau),u(\tau),p(\tau))$ coincide with $(X(\tau),U(\tau),P(\tau))$ at the maximum interval $[t_0,t_0+\delta)$ where $\delta<t-t_0$. By (L3), there exists a finite constant $K$ such that for any $\tau\in [t_0,t_0+\delta)$,
$|\dot{\bar{\gm}}(\tau)|<K$.
Hence, we have
\[|\bar{\gm}(t_0+\delta)-\bar{\gm}(t_0+\delta/2)|\leq \int_{t_0+\delta/2}^{t_0+\delta}|\dot{\bar{\gm}}(\tau)|d\tau<K\delta/2.\]
Let
\[F(\tau)=|\bar{\gm}(\tau)-\bar{\gm}(t_0+\delta/2)|-K(\tau-(t_0+\delta/2)).\]
Since $F(t_0+\delta)<0$, it follows from the continuity of $F$ that for $s>t_0+\delta$ and close to $t_0+\delta$, we also have
$F(s)<0$. That is
\[|\bar{\gm}(s)-\bar{\gm}(t_0+\delta/2)|<K(s-(t_0+\delta/2)).\]
Hence, $(s,\bar{\gm}(s))$ is contained in a cone-like region. Repeat the arguments in Step (a) and Step (b), it follows  that for any $\tau\in [t_0,s)$,
 \[(\bar{\gm}(\tau),u(\tau),p(\tau))=(X(\tau),U(\tau),P(\tau)),\] which is in contradiction with the assumption.
 Therefore, for any $\tau \in (0,t)$, there holds
\[(\bar{\gm}(\tau),u(\tau),p(\tau))=(X(\tau),U(\tau),P(\tau)).\]
 We complete the proof of Lemma \ref{chc}.\End

So far, we complete the proof of Theorem \ref{two} under the assumptions (A1)-(A2) and (L1)-(L3). \End

%In the following, we denote
%\begin{equation}
%h(x_0,0):=\lim_{t\rightarrow 0^+}U(t)=u_0.
%\end{equation}

%\begin{Remark}
%In order to avoid the tendinous argument of the classical characteristics method, we use Malgrange preparation theorem to obtain the classical solution $S(x,t)$ in a cone-like region rather than a cuboid-like region, which is sufficient for the deduction.  But we have to restrict ourselves to the $C^\infty$ systems due to that theorem. It is not difficult to check that $C^r$ ($r\geq 2$) systems can be  handled similarly  by virtue of the classical characteristics method and Theorem \ref{two1} can be obtained more directly.
%\end{Remark}
\subsection{Solution and action function}
In this step,  we will prove Theorem \ref{two2} under the assumptions (A1)-(A2) and (L1)-(L3).

In order to prove Theorem \ref{two2}, we need the following lemma:
\begin{Lemma}\label{baru}
There exists a solution of (\ref{hjech}) denoted by $(\bar{X}(t),\bar{U}(t),\bar{P}(t))$ such that
\begin{equation}\label{ifff}
\bar{U}(t)=\inf\left\{U(t):(X(s),U(s),P(s))\in \mathcal{S}_{x_0,u_0}^x(t)\right\}.
\end{equation}
\end{Lemma}
\Proof
The proof is divided into three claims. The first one shows the infimum of (\ref{ifff}) is finite.

\noindent\textbf{Claim A:} There exists a finite $\bar{U}_t\in \R$ such that
\begin{equation}
\bar{U}_t=\inf_{\mathcal{S}_{x_0,u_0}^x(t)}U(t).
\end{equation}

\noindent\textbf{Proof of Claim A:}
In terms of the contact Hamiltonian equation (\ref{hjech}), it follows from the Legendre transformation that
\begin{equation}\label{U1}
\dot{U}(t)=L(X(t),U(t),\dot{X}(t)).
\end{equation}
By the assumption (A1),  for any $s\in [0,t]$,
\[L(X(s),U(s),\dot{X}(s))\geq L(X(s),0,\dot{X}(s))-\lambda|U(s)|.\]

Let us estimate the lower bound of $U(s)$ for $s\in [0,t]$. Consider the set
\[\mathcal{Z}:=\{\tau\ |\ U(\tau)<0, \ \tau\in [0,t]\}.\]
Without loss of generality, we assume $\mathcal{Z}\neq\emptyset$, otherwise $U(s)\geq 0$ for any $s\in [0,t]$  which verifies the claim. For any $s\in \mathcal{Z}$, we have two cases as follows.

Case One: if $\sup_{\tau\in [0,s)}U(\tau)\geq 0$, then there exits $s_0\in [0,s)$ such that for any $\tau\in [s_0,s]$, $U(\tau)\leq 0$ and $U(s_0)=0$. By integration, we have
\begin{equation}
U(s)=U(s)-U(s_0)\geq \int_{s_0}^sL(X(\tau),0,\dot{X}(\tau))d\tau+\lambda\int_{s_0}^s U(\tau)d\tau.
\end{equation}
Let $V(s):=-U(s)$ for $s\in [0,t]$. There holds
\begin{equation}
V(s)\leq -\int_{s_0}^sL(X(\tau),0,\dot{X}(\tau))d\tau+\lambda\int_{s_0}^s V(\tau)d\tau.
\end{equation}
By (L2), $L(X(\tau),0,\dot{X}(\tau))$ is bounded from below. Hence,
\begin{equation}
V(s)\leq C+\lambda\int_{s_0}^s V(\tau)d\tau,
\end{equation}
which together with Gronwall's inequality implies $V(s)\leq Ce^{\lambda (s-s_0)}$. Hence, we have
\begin{equation}\label{ust}
U(s)\geq-Ce^{\lambda (s-s_0)}\geq -Ce^{\lambda t }.
\end{equation}

Case Two: if $\sup_{\tau\in [0,s)}U(\tau)< 0$, a similar argument yields
\[U(s)\geq-Ce^{\lambda s}\geq -Ce^{\lambda t}.\]
Therefore, $\inf_{\mathcal{S}_{x_0,u_0}^x}U(t)$ exists, which is denoted by $\bar{U}_t$. We complete the proof of Claim A.\End

  Then, one can find a sequence $(X_n(t),U_n(t),\dot{X}_n(t))$ such that (extracting a subsequence if necessary)
$U_n(t)\rightarrow \bar{U}_t$  as $n\rightarrow\infty$.

\noindent\textbf{Claim B:} For any $s\in [0,t]$ and $n$ large enough, we have
\begin{equation}\label{unsc}
|U_n(s)|\leq C,
\end{equation}
where $C$ is a constant independent of $n$.

\noindent\textbf{Proof of Claim B:}
 By (\ref{ust}), $U_n(s)$ is bounded from below. In the following, we prove $U_n(s)$ is bounded from above. We denote $K_n:=|U_n(t)|$.
For $n$ large enough, we have
\begin{equation}\label{U2}
K_n\leq |\bar{U}_t|+1.
\end{equation}
 We prove $U_n(s)$ is bounded from above for any $s\in [0,t]$ and $n$ large enough. Let
 \[\delta=\frac{1}{2\lambda}.\]
 The proof is divided into two cases.

Case One: $t\leq \delta$.  We prove there exists a suitable large constant $l>0$ such that for any $s\in [0,t]$, $U_n(s)<lK_n\leq l|\bar{U}_t|+l$.

By contradiction, we assume that there exists $s_0$ such that $U_n(s_0)\geq lK_n$. Hence, one can find a time interval $[s_1,s_2]\subset[0,t]$ such that $U_n(s_1)=lK_n$, $U_n(s_2)=K_n$ and $K_n\leq U_n(s)\leq lK_n$ for $s\in [s_1,s_2]$. Moreover, we have
\begin{align*}
U_n(s_2)&=U_n(s_1)+\int_{s_1}^{s_2}L(X_n(\tau),U_n(\tau),\dot{X}_n(\tau))d\tau,\\
&\geq U_n(s_1)-C(s_2-s_1)-\lambda\int_{s_1}^{s_2}U_n(\tau)d\tau,
\end{align*}
where $C$ denotes the lower bound of $L(x,0,\dot{x})$. Combining with $t\leq \delta=\frac{1}{2\lambda}$, it yields
\[U_n(s_2)\geq \frac{1}{2}lK_n-C\delta,\]
which is in contradiction with the assumption $U_n(s_2)=K_n$ for $l$ suitable large. Hence, for any $s\in [0,t]$, $U_n(s)<lK_n\leq l|\bar{U}_t|+l$.

Case Two: $t>\delta$. In this case, it is easy to see that there exits $m>0$ such that the time interval is separated into $m$ subintervals and the length of each one is not greater than $\delta$. More precisely, one can choose a partition:
\[\left\{\tau_i : \tau_{i+1}-\tau_i\leq \frac{\delta}{2}\right\},\]
where $i=0,1,\ldots,2m$ and $\tau_0=0, \tau_{2m}=t$. Hence, for any $\bar{s}\in (0,t)$, there exists $\bar{i}\geq 1$ such that $\bar{s}\in (\tau_{\bar{i}-1},\tau_{\bar{i}+1}]$. By an argument similar to Case One, we have
\[U_n(\bar{s})<lU_n(\tau_{\bar{i}+1})<l^2U_n(\tau_{\bar{i}+2}).\]
Repeating the deduction by $2m-\bar{i}$ times, it follows that
\[U_n(\bar{s})<l^{2m-\bar{i}}U_n(\tau_{2m})=l^{2m-\bar{i}}U_n(t)\leq l^{2m}K_n.\]
We denote $\bar{l}:=l^{2m}$. It follows that there exists $\bar{l}>0$ suitable large such that for any $s\in [0,t]$,
$U_n(s)<\bar{l}K_n\leq \bar{l}|\bar{U}_t|+\bar{l}$.
Therefore, for both of the cases, one can find a constant $C$ independent of $n$ such that
$|U_n(s)|\leq C$. We complete the proof of Claim B.\End

\noindent\textbf{Claim C:} There holds
\begin{equation}
|\dot{X}_n(0)|\leq C,
\end{equation}
where $C$ is a constant independent of $n$.

\noindent\textbf{Proof of Claim C:} By contradiction, we assume that there exists a subsequence $|\dot{X}_{n}(0)|$ such that
\[|\dot{X}_{n}(0)|\geq n.\]
 According to (\ref{U1}), it follows from (A1) that
\begin{align*}
U_n(t)-U_n(0)&=\int_{0}^{t}L(X_n(s),U_n(s),\dot{X}_n(s))ds\\
&\geq \int_{0}^{t}L(X_n(s),0,\dot{X}_n(s))ds-\lambda\int_{0}^{t} |U(s)|ds,
\end{align*}
which together with (\ref{unsc}) and (L2) implies
\[\int_{0}^{t}|\dot{X}_n(s)|ds \leq C.\]
Since $\dot{X}_n(s)$ is continuous with respect to $s$, by the mean value theorem for integration,
for any $n\in\N$, one can find a  $s_0^n\in [0,t]$ satisfying $|\dot{X}_n(s_0^n)|\leq C$.  Extracting a subsequence if necessary, it follows from (\ref{unsc}) and the compactness of $M$ that
\[s_0^n\rightarrow s_0,\quad (X_n(s_0^n),U_n(s_0^n),\dot{X}_n(s_0^n))\rightarrow (\bar{X},\bar{U},\bar{V}), \quad n\rightarrow\infty.\]
Let $\Phi_s$ be the flow generated by $L(x,u,\dot{x})$. Then it follows from the completeness of the flow that for any $s\in [0,t]$, $\Phi_{s}(\bar{X},\bar{U},\bar{V})$ is well defined.  We consider $(\bar{X},\bar{U},\bar{V})$ as the initial condition of $\Phi_s$. Based on the continuous dependence of solutions of ODEs on initial conditions, it follows that for $n$ large enough,
\[\text{dist}\left((X_n(0),U_n(0),\dot{X}_n(0)),\Phi_{-s_0}(\bar{X},\bar{U},\bar{V}))\right)\leq 1,\]
where dist$(\cdot,\cdot)$ denotes the distance induced by the Riemannian metric on $TM\times\R$.
Hence, $|\dot{X}_n(0)|$ has a bound independent of $n$. We complete the proof of Claim C.\End

\noindent\textbf{Proof of Lemma \ref{baru}:}
By Claim C, for $n$ large enough,
\begin{equation}
|\dot{X}_n(0)|\leq C.
\end{equation}
It follows that extracting a subsequence if necessary,  as $n\rightarrow\infty$
\begin{equation}
\dot{X}_n(0)\rightarrow \bar{v}, \bar{U}_n(t)\rightarrow \bar{U}_t.
\end{equation}
In addition, we have
\[X_n(0)=x_0, U_n(0)=u_0, X_n(t)=x.\]
By virtue of continuous dependence of solutions of ODEs on initial conditions, it follows that there exists an orbit generated by $L$ (via Legendre transformation) $(\bar{X}(t),\bar{U}(t),\dot{\bar{X}}(t))$  such that
\begin{equation}
\begin{array}{lll}
\bar{X}(0)=x_0, \bar{U}(0)=u_0, \dot{\bar{X}}(0)=\bar{v}, \bar{X}(t)=x, \bar{U}(t)=\bar{U}_t.\\
\end{array}
\end{equation}Therefore, there exists $(\bar{X}(t),\bar{U}(t),\bar{P}(t))$ such that
\begin{equation}\label{inf}
\bar{U}(t)=\inf\left\{U(t):(X(s),U(s),P(s))\in \mathcal{S}_{x_0,u_0}^x(t)\right\}.
\end{equation}
This completes the proof of Lemma \ref{baru}.
\End

\noindent\textbf{Proof of Theorem \ref{two2}:}
In the following, we will prove
\begin{equation}
h_{x_0,u_0}(x,t)=\bar{U}(t).
\end{equation}

By virtue of (\ref{U1}), we have
\begin{equation}
\bar{U}(t)=u_0+\int_0^tL(\bar{X}(\tau),\bar{U}(\tau),\dot{\bar{X}}(\tau))d\tau.
\end{equation} By Lemma \ref{exist} and Lemma \ref{chc}, it follows that there exists an orbit generated by $L$ (via Legendre transformation) $(\hat{X}(t),\hat{U}(t),\dot{\hat{X}}(t))$ such that  $\hat{U}(t)=h_{x_0,u_0}(x,t)$ and
\begin{equation}\label{hxu}
h_{x_0,u_0}(x,t)=u_0+\int_0^tL(\hat{X}(\tau),h_{x_0,u_0}(\hat{X}(\tau),\tau),\dot{\hat{X}}(\tau))d\tau.
\end{equation}
Hence, $\bar{U}(t)\leq h_{x_0,u_0}(x,t)$. It suffices to prove $\bar{U}(t)\geq h_{x_0,u_0}(x,t)$.
By contradiction, we assume   $\bar{U}(t)<h_{x_0,u_0}(x,t)$.
 By (\ref{U1}), we have
\begin{align*}
\bar{U}(t)=u_0+\int_0^{t}L(\bar{X}(\tau),\bar{U}(\tau),\dot{\bar{X}}(\tau))d\tau.
\end{align*}
\begin{align*}
h_{x_0,u_0}(x,t)\leq u_0+
\int_0^{t}L(\bar{X}(\tau),h_{x_0,u_0}(\bar{X}(\tau),\tau),\dot{\bar{X}}(\tau))d\tau.
\end{align*}
For any $\sigma\in [0,t]$, we denote $\bar{H}(\sigma):=h_{x_0,u_0}(\bar{X}(\sigma),\sigma)$. In particular, we have $\bar{H}(t):=h_{x_0,u_0}(\bar{X}(t),t)$.  Let
 \begin{equation}
F(\sigma):=\bar{H}(\sigma)-\bar{U}(\sigma),
\end{equation}where $\sigma\in [0,t]$. It is easy to see that $\bar{U}(0)=u_0=\bar{H}(0)$. Then we have  $F(0)=0$. The assumption $\bar{U}(t)<h_{x_0,u_0}(x,t)$ implies $F(t)>0$. Hence, there exists $\sigma_0\in [0,t)$ such that $F(\sigma_0)=0$ and  $F(\sigma)> 0$ for $\sigma> \sigma_0$.
Moreover, for any $\tau\in (\sigma_0,t]$, we have
  \begin{equation}\label{s77}
\bar{U}(\tau)=\bar{U}(\sigma_0)+\int_{\sigma_0}^{\tau}L(\bar{X}(\sigma),\bar{U}(\sigma),\dot{\bar{X}}(\sigma))d\sigma.
\end{equation}
Let $\tilde{X}$ be a calibrated curve of $h_{x_0,u_0}$ with $\tilde{X}(0)=x_0$, $\tilde{X}(\sigma_0)=\bar{X}(\sigma_0)$. For $\sigma\in [\sigma_0,\tau]$, we construct a family of $X_\tau:[0,\tau]\rightarrow M$ as follows:
\begin{equation}
X_\tau(\sigma)=\left\{\begin{array}{ll}
\hspace{-0.4em}\tilde{X}(\sigma),& \sigma\in [0,\sigma_0],\\
\hspace{-0.4em}\bar{X}(\sigma),&\sigma\in (\sigma_0,\tau].\\
\end{array}\right.
\end{equation}
Moreover, for any $\tau\in (\sigma_0,t]$, we have
  \begin{equation}
\bar{H}(\tau)\leq \bar{H}(\sigma_0)+\int_{\sigma_0}^{\tau}L(\bar{X}(\sigma),\bar{H}(\sigma),\dot{\bar{X}}(\sigma))d\sigma.
\end{equation}
Since $\bar{H}(\sigma_0)-\bar{U}(\sigma_0)=F(\sigma_0)=0$, a direct calculation implies
  \begin{equation}
\bar{H}(\tau)-\bar{U}(\tau)\leq \int_{\sigma_0}^\tau\lambda(\bar{H}(\sigma)-\bar{U}(\sigma))d\sigma.
\end{equation}
Hence, we have
  \begin{equation}
F(\tau)\leq \int_{\sigma_0}^\tau\lambda F(\sigma)d\sigma,
\end{equation}
which is
 in contradiction with Lemma \ref{gron}. Hence, we have
 \begin{equation}\label{s22uh}
\bar{U}(t)\geq h_{x_0,u_0}(x,t).
\end{equation} It contradicts the assumption $\bar{U}(t)<h_{x_0,u_0}(x,t)$.
 This finishes the proof of Theorem \ref{two2} under the assumptions (A1)-(A2) and (L1)-(L3).
\End

\section{\sc Back to the original Lagrangian}

 In the following, we  focus on the relaxation of the assumptions from the uniform Lipschitzity and uniform boundedness to the Osgood growth and prove the theorems for the original Lagrangian under (L1)-(L4).
 The proof will be divided into two steps. In the first step,  we will prove a priori estimate of the  orbit generated by the modified Lagrangian. In the second step, we will introduce the action function of the original Lagrangian. In the last step, the proofs of the theorems for the original Lagrangian   will be completed by an argument on a limiting passage.

 Let $H_R$ denote the Hamiltonian associated to $L_R$ via the Legendre transformation.
 Based on the arguments in Section 3,  we have
 \begin{Lemma}\label{inter}
For given $x_0\in M$, $u_0\in\R$ and $T>0$, there exists a  unique continuous function $h_{x_0,u_0,R}(x,t)$ defined on $M\times (0,T]$ satisfying
\begin{equation}
h_{x_0,u_0,R}(x,t)=u_0+\inf_{\substack{\gm(t)=x \\  \gm(0)=x_0} }\int_0^tL_R(\gm(\tau),h_{x_0,u_0,R}(\gm(\tau),\tau),\dot{\gm}(\tau))d\tau,
\end{equation}where the infimum is taken among the absolutely continuous curves $\gm:[0,t]\rightarrow M$. In particular, the infimum is attained at  a $C^1$ curve denoted by $\bar{\gm}$. Let
\[x(t):=\bar{\gm}(t),\quad u(t):=h_{x_0,u_0,R}(\bar{\gm}(t),t),\quad p(t)=\frac{\partial L_R}{\partial \dot{x}}(\bar{\gm}(t),h_{x_0,u_0,R}(\bar{\gm}(t),t),\dot{\bar{\gm}}(t)),\]
then $(x(t),u(t),p(t))$ satisfies the contact Hamiltonian equation generated by $H_R$.
 Moreover, let $\mathcal{S}_{x_0,u_0,R}^x(t)$ denote the set of the solution  of (\ref{hjech}) generated by $H_R$ satisfying $X(0)=x_0$, $X(t)=x$ and $U(0)=u_0$, then we have
\begin{equation}
h_{x_0,u_0,R}(x,t)=\inf\left\{U(t):(X(s),U(s),P(s))\in \mathcal{S}_{x_0,u_0,R}^x(t)\right\}.
\end{equation}
 \end{Lemma}

We will omit the subscripts $x_0, u_0$ and denote $h_{R}(x,t):=h_{x_0,u_0,R}(x,t)$.

\subsection{A priori compactness}
In this step, we will prove a priori estimate of $|h_R(\gm_R(s),s)|$ and $|\dot{\gm}_R(s)|$, where $\gm_R:[0,t]\rightarrow M$ is a calibrated curve of $h_R$ satisfying  $\gm_R(0)=x_0$ and $\gm_R(t)=x$.
\begin{Lemma}\label{flowcomplete}
Let $(X_R(s),U_R(s),P_R(s))$ denote the solution  of (\ref{hjech}) generated by $H_R$ satisfying $X_R(s)=x_0$, $X_R(t)=x$ and $U_R(s)=u_0$. Let
\[ \dot{X}_R(t)=\frac{\partial H_R}{\partial p}(X_R(t),U_R(t),P_R(t)).,\]
If there exist $s_0\in [0,t]$ and a positive constant $C$ independent of $R$ such that
\begin{equation}\label{xrsoo}
|U_R(s_0)|\leq C,\quad |\dot{X}_R(s_0)|\leq C,
\end{equation}
then for any $s\in [0,t]$, we have
\[(X_R(s),U_R(s),\dot{X}_R(s))\in \mathcal{K},\]
where $\mathcal{K}$ is a compact set independent of $R$.
\end{Lemma}
\Proof
By virtue of the completeness of the flow (L3),  it suffices to prove that for any $R$,
$|\dot{X}_R(0)|$ has a bound independent of $R$ since $X_R(s)=x_0$ and $U_R(s)=u_0$.

By contradiction, we assume that there exists a subsequence $R_n$ such that
\begin{equation}
|\dot{X}_{R_n}(0)|\geq n.
\end{equation}
For the simplicity of notations, we still use $|\dot{X}_{R}(0)|$ to denote the subsequence.

 Extracting a subsequence if necessary, it follows from (\ref{xrsoo}) that
 \[|U_R(s_0^n)|\leq C,\quad |\dot{X}_R(s_0^n)|\leq C,\]
which together with the compactness of $M$ yields
\[s_0^n\rightarrow s_0,\quad (X_R(s_0^n),U_R(s_0^n),\dot{X}_R(s_0^n))\rightarrow (\bar{X},\bar{U},\bar{V}), \quad n\rightarrow\infty.\]
Let $\Phi^R_s$ be the flow generated by $L_R(x,u,\dot{x})$. It follows from the completeness of the flow that for any $s\in [0,t]$, $\Phi^R_{s}(\bar{X},\bar{U},\bar{V})$ is well defined.  We consider $(\bar{X},\bar{U},\bar{V})$ as the initial condition of $\Phi^R_s$. Based on the continuous dependence of solutions of ODEs on initial conditions, it follows that for $R$ large enough,
\[\text{dist}\left((X_R(0),U_R(0),\dot{X}_R(0)),\Phi_{-s_0}(\bar{X},\bar{U},\bar{V}))\right)\leq 1.\]
Hence, $|\dot{X}_R(0)|$ has a bound independent of $R$.
\End

\begin{Lemma}\label{aprioriset}
Let $(X_R(s),U_R(s),\dot{X}_R(s))$ denote the orbit generated by $L_R$ satisfying $X_R(0)=x_0$, $X_R(t)=x$ and $U_R(0)=u_0$. For each $C_1>0$, if  $U_R(t)\leq C_1$,
then there exists a positive constant $C_2:=C_2(C_1)$ independent of $R$ such that  for any $s\in [0,t]$,
\begin{equation}
|U_R(s)|\leq C_2.
\end{equation}
\end{Lemma}
\Proof
On one hand, we prove $U_R(s)$ is
 bounded from below.  By contradiction, we assume that extracting a subsequence if necessary, there exists  subsequences $R_n$ and $s_n$ such that
\begin{equation}\label{81}
U_{R_n}(s_n)<-n.
\end{equation}
 Since $U_R(0)=u_0$ for any $R$, then it follows from (\ref{81}) that for $n$ large enough, there exists a time interval $[t_1,t_2]$ such that
\begin{equation}\label{uooo}
\left\{\begin{array}{ll}
\hspace{-0.4em}U_{R_n}(t_1)=u_0,& t_1\in [0,t),\\
\hspace{-0.4em}U_{R_n}(t_2)=u_0-1,& t_2\in (0,t],\\
\hspace{-0.4em}U_{R_n}(s)\in [u_0-1,u_0], &\text{for}\  s\in[t_1,t_2],\\
\hspace{-0.4em}U_{R_n}(s)\geq u_0-1, &\text{for}\  s\in[0,t_2].\\
\end{array}\right.
\end{equation}
 Indeed, by the continuity of $U_{R_n}(s)$ with respect to $s$, it suffices to take $n>1-u_0$.

By Lemma \ref{inter}, we have
\[\dot{U}_{R_n}(s)=L_{R_n}(X_{R_n}(s),U_{R_n}(s),\dot{X}_{R_n}(s)).\]
By (\ref{uooo}), there exists $s_0\in [t_1,t_2]$ such that $\dot{U}_{R_n}(s_0)<0$.
Since $U_{R_n}(s)\in [u_0-1,u_0]$ for $s\in [t_1,t_2]$, it follows from (\ref{import}) that there exists a constant $B$ independent of $R_n$ such that
\[\dot{U}_{R_n}(s_0)\geq |\dot{X}_{R_n}(s_0)|-B.\]
Hence, we have
\[|\dot{X}_{R_n}(s_0)|<B.\]
Moreover, the assumptions of Lemma \ref{flowcomplete} are satisfied. It follows that for any $s\in [0,t]$,
\[(X_{R_n}(s),U_{R_n}(s),\dot{X}_{R_n}(s))\in \mathcal{K},\]
where $\mathcal{K}$ is independent of $R_n$. In particular, it is in contradiction with (\ref{81}). Therefore, we have
\begin{equation}
U_{R}(s)\geq C.
\end{equation}

On the other hand, we prove $U_{R}(s)$ is  bounded from above.
We denote $u_t:=U_R(t)$. By contradiction, we assume that there exists  subsequences $R_n$ and $s_n$ such that
\begin{equation}\label{84}
U_{R_n}(s_n)>n.
\end{equation}
It follows from (\ref{84}) and the continuity of $U_{R_n}(s)$ that for $n$ large enough, there exists a time interval $[t_1,t_2]$ such that
 \begin{equation}\label{uooo1}
\left\{\begin{array}{ll}
\hspace{-0.4em}U_{R_n}(t_1)=u_t+1,& t_1\in [0,t),\\
\hspace{-0.4em}U_{R_n}(t_2)=u_t,& t_2\in (0,t],\\
\hspace{-0.4em}U_{R_n}(s)\in [u_t,u_t+1], &\text{for}\  s\in[t_1,t_2],\\
\hspace{-0.4em}U_{R_n}(s)\leq u_t+1, &\text{for}\  s\in[t_1,t].\\
\end{array}\right.
\end{equation}
 By Lemma \ref{inter}, we have
\[\dot{U}_{R_n}(s)=L_{R_n}(X_{R_n}(s),U_{R_n}(s),\dot{X}_{R_n}(s)).\]
By (\ref{uooo1}), there exists $s_0\in [t_1,t_2]$ such that $\dot{U}_{R_n}(s_0)<0$.
Since $U_{R_n}(s)\in [u_t,u_t+1]$ for $s\in [t_1,t_2]$, it follows from (\ref{import}) that there exists a constant $B$ independent of $R_n$ such that
\[\dot{U}_{R_n}(s_0)\geq |\dot{X}_{R_n}(s_0)|-B.\]
Hence, we have
\[|\dot{X}_{R_n}(s_0)|<B.\]
Moreover, the assumptions of Lemma \ref{flowcomplete} are satisfied.
It follows that for any $s\in [0,t]$,
\[(X_{R_n}(s),U_{R_n}(s),\dot{X}_{R_n}(s))\in \mathcal{K},\]
where $\mathcal{K}$ is independent of $R_n$.  In particular, it is in contradiction with (\ref{84}). Therefore, we have
\begin{equation}
U_R(s)\leq C.
\end{equation}
 We complete the proof of Lemma \ref{aprioriset}.
\End

Based on Lemma \ref{aprioriset}, we can obtain a priori estimate of  $|h_R(\gm_R(s),s)|$.
\begin{Lemma}[Pointwise Boundedness]\label{A}
For a given $(x,t)\in M\times (0,T]$, let $\gm_R(s):[0,t]\rightarrow M$ be a calibrated curve of $h_R$ satisfying $\gm_R(0)=x_0$ and $\gm_R(t)=x$, then  there exists a positive constant $C$ independent of $R$ such that for any $s\in [0,t]$,
\begin{equation}\label{astar}
|h_R(\gm_R(s),s)|\leq C.
\end{equation}
\end{Lemma}
\Proof
In terms of Lemma \ref{inter}, there holds for any $s\in [0,t]$,
\[(\gm_R(s),h_R(\gm_R(s),s),\dot{\gm}_R(s))=(X_R(s),U_R(s),\dot{X}_R(s)).\]
By Lemma \ref{aprioriset}, it suffices to show that for a given $(x,t)\in M\times (0,T]$,
  $h_R(x,t)$ is  bounded from above.

 Without loss of generality, we assume $h_R(x,t)>0$, otherwise, we obtain the upper bound of $h_R(x,t)$ directly. Let $\bar{\gm}(s):[0,t]\rightarrow M$ be a straight line satisfying $\bar{\gm}(0)=x_0$ and $\bar{\gm}(t)=x$.
 Since $h_R(\bar{\gm}(t),t)=h_R(x,t)>0$,  then we have the following two cases.

Case One:
  there exists $s_0\in (0,t)$ such that $h_R(\bar{\gm}(s_0),s_0)=0$ and $h_R(\bar{\gm}(s),s)\geq 0$ for any $s\in [s_0,t]$.

We denote $h_R(s):=h_R(\bar{\gm}(s),s)$ for the sake of the simplicity. Since $\bar{\gm}$ is a straight line, then it follows from the compactness that for $s\in [0,t]$, $(\bar{\gm}(s),\dot{\bar{\gm}}(s))$ is contained in a compact set independent of $R$ denoted by $K$.
Moreover,  we have
\begin{align*}
h_R(t)\leq h_R(s_0)+\int_{s_0}^{t}L_R(\bar{\gm}(\tau),h_R(\tau),\dot{\bar{\gm}}(\tau))d\tau
\leq \int_{s_0}^{t}f_K(h_R(\tau))d\tau,
\end{align*}where the second inequality is from the Osgood growth assumption (L4).
Let $g_R(\tau)$ be a function defined on $[0,t-s_0]$ and satisfy
\begin{equation}
\begin{cases}
\dot{g}_R(\tau)= f_K(g_R(\tau)),\\
g_R(0)=h_R(s_0)=0.
\end{cases}
\end{equation}
Hence, we have
\begin{equation}
\int_0^{g_R}\frac{1}{f_K(g_R)}dg_R=\int_0^sd\tau.
\end{equation}
In particular, we have
\begin{equation}\label{aosgo}
\int_0^{g_R(t-s_0)}\frac{1}{f_K(g_R)}dg_R=t-s_0,
\end{equation}which together with (L4) yields $g_R(t-s_0)$ has a upper bound independent of $R$. Otherwise, the left side of (\ref{aosgo}) diverges, which contradicts the right side (\ref{aosgo}). By the comparison theorem of ODEs (see \cite{HS} for instance), we have $h_R(t)\leq g_R(t-s_0)$.  Hence, $h_R(\gm_R(t),t)$ has a upper bound independent of $R$.

 Case Two: for any $s\in (0,t)$, $h_R(\bar{\gm}(s),s)>0$.

 By a similar argument as the previous case, we have $h_R(t)\leq g_R(t)$, where $g_R$ has a upper bound independent of $R$ and satisfies
 \begin{equation}
\begin{cases}
\dot{g}_R(\tau)= f_K(g_R(\tau)),\\
g_R(0)=h_R(0)=u_0.
\end{cases}
\end{equation}

This finishes the proof of Lemma \ref{A}.
\End

%\begin{Remark}
%The Osgood growth assumption (H4) ((L4) via Legendre transformation) is only used to obtain the upper bound of  $h_R(\gm_R(t),t)$. It seems to be inevitable that this assumption can be improved.
%\end{Remark}

For a given $k>0$, let
\[\Sigma_k:=\{(x,t)\in M\times[0,T]\ |\ d(x,x_0)\leq kt\},\]
where $d(\cdot,\cdot)$ denotes the distance between $x$ and $x_0$ induced by the Riemannian metric on $M$.
Based on a similar argument as the one in the proof of Lemma \ref{A}, we have the following lemma.
\begin{Lemma}[Uniform Boundedness]\label{Auuu}
For any $(x,t)\in \Sigma_k$, let $\gm_R(s):[0,t]\rightarrow M$ be a calibrated curve of $h_R$ satisfying $\gm_R(0)=x_0$ and $\gm_R(t)=x$, then there exists a positive constant $C$ independent of $R$ such that   for any $s\in [0,t]$,
\begin{equation}\label{astar}
|h_R(\gm_R(s),s)|\leq C.
\end{equation}
\end{Lemma}

Based on Lemma \ref{flowcomplete}, Lemma \ref{A} and Lemma \ref{Auuu}, we obtain a priori compactness.
\begin{Lemma}\label{xr}
Let $\gm_R:[0,t]\rightarrow M$ be a calibrated curve of $h_R$ satisfying $\gm_R(0)=x_0$ and $\gm_R(t)=x$.
\begin{itemize}
\item [(i)] For a given $(x,t)\in M\times(0,T]$, there exists a compact set $\mathcal{K}$  independent of $R$ such that for any $R$ and $s\in [0,t]$, there holds
\[(\gm_R(s),h_R(\gm_R(s),s),\dot{\gm}_R(s))\in \mathcal{K},\]
where $\mathcal{K}$ is a compact set independent of $R$.
\item [(ii)] For any $(x,t)\in \Sigma_k$, there exists a compact set $\mathcal{K}'$  independent of $R$ such that for any $R$ and $s\in [0,t]$, there holds
\[(\gm_R(s),h_R(\gm_R(s),s),\dot{\gm}_R(s))\in \mathcal{K}',\]
where $\mathcal{K}'$ is a compact set independent of $R$.
\end{itemize}
\end{Lemma}
\Proof
It suffices to prove (i). (ii) can be obtained by a similar argument as the one in Lemma \ref{Auuu}.
In terms of Lemma \ref{inter}, there holds for any $s\in [0,t]$,
\[(\gm_R(s),h_R(\gm_R(s),s),\dot{\gm}_R(s))=(X_R(s),U_R(s),\dot{X}_R(s)).\]
Clearly, we have
\[\dot{U}_{R}(s)=L_{R}(X_{R}(s),U_{R}(s),\dot{X}_{R}(s)).\]
By Lemma \ref{A}, we have
\begin{equation}\label{urscc}
|U_R(s)|\leq C,
\end{equation}
which implies that
\begin{equation}\label{822aa}
\int_{0}^{t}L_{R}(X_{R}(s),U_{R}(s),\dot{X}_{R}(s))ds=U_{R}(t)-U_{R}(0)\leq C-u_0.
\end{equation}
By (\ref{import}), we have
\[\int_{0}^{t}|\dot{X}_{R}(s)|-Bds\leq C-u_0.\]
By the mean value theorem for integration, it follows that there exists $s_0\in [0,t]$ such that for any $R$,
\begin{equation}\label{xrscb}
|\dot{X}_{R}(s_0)|\leq \frac{C-u_0}{t}+B.
\end{equation}
Hence, Lemma \ref{xr} follows from Lemma \ref{flowcomplete} and Lemma \ref{A}.
\End

\subsection{Action function}

In this step, we will introduce the action function of the original Lagrangian. For a given $R>0$, we denote
\[\mathcal{D}_R=\{(x,u,\dot{x})\ :\ x\in M, |u|\leq R, |\dot{x}|\leq R\}.\]
By Lemma \ref{xr}, for a given $(x,t)\in M\times (0,T]$, there exists $R_0$ such that for any $R>R_0$, we have
\begin{equation}\label{R00}
(\gm_R(s),h_R(\gm_R(s),s),\dot{\gm}_R(s))\in \mathcal{K}\subset \mathcal{D}_{R_0},
\end{equation}where $\gm_R:[0,t]\rightarrow M$ is a calibrated curve of $h_R$ satisfying $\gm_R(0)=x_0$ and $\gm_R(t)=x$.
As a preliminary, we have the following lemma:

\begin{Lemma}[Pointwise Invariance]\label{heh} For a given $(x,t)\in M\times (0,T]$ and for any $R_1, R_2>R_0$, there holds
\begin{equation}
h_{R_1}(x,t)=h_{R_2}(x,t)=h_{R_0}(x,t),
\end{equation}where $R_0$ is determined by (\ref{R00}).
\end{Lemma}
\Proof
 Let $(X_R(s),U_R(s),P_R(s))$ denote a solution  of (\ref{hjech}) generated by $H_R$, where $H_R$ denotes the Legendre transformation of $L_R$.  By Lemma \ref{inter}, we have
\begin{equation}\label{xxrr}
h_{R}(x,t)=\inf\left\{U_R(t):(X_R(s),U_R(s),P_R(s))\in \mathcal{S}_R\right\},
\end{equation}
 where we omit the subscripts $x_0,u_0,x,t$ of $\mathcal{S}_{x_0,u_0,R}$. According to Lemma \ref{baru}, there exists a solution  $(\bar{X}_{R}(t),\bar{U}_{R}(t),\bar{P}_{R}(t))$ of (\ref{hjech}) such that
\begin{equation}
\bar{U}_{R}(t)=\inf\left\{U_R(t):(X_R(s),U_R(s),P_R(s))\in \mathcal{S}_{R}\right\},
\end{equation}which together with (\ref{xxrr}) implies that for $(x,t)\in M\times (0,T]$,
\[h_{R}(x,t)=\bar{U}_{R}(t).\]
It suffices to prove for any $R_1, R_2>R_0$, $\bar{U}_{R_1}(t)=\bar{U}_{R_2}(t)$. In fact, it follows from  (\ref{R00}) that for any $R>R_0$,
\[(\bar{X}_R(s),\bar{U}_R(s),\bar{P}_R(s))\in \mathcal{D}^*_{R_0},\]
where $\mathcal{D}^*_{R_0}$ denotes the dual set of $\mathcal{D}_{R_0}$ associated to the Legendre transformation.
Hence, for any $R_1, R_2>R_0$,  we have
\begin{align*}
\bar{U}_{R_1}(t)&=\inf\left\{U_R(t):(X_R(s),U_R(s),P_R(s))\in \mathcal{S}_{R_1}\right\},\\
&=\inf\left\{U_R(t):(X_R(s),U_R(s),P_R(s))\in \mathcal{S}_{R_1}\cap\mathcal{D}^*_{R_0}\right\},\\
&=\inf\left\{U_R(t):(X_R(s),U_R(s),P_R(s))\in \mathcal{S}_{R_2}\cap\mathcal{D}^*_{R_0}\right\},\\
&=\bar{U}_{R_2}(t)=\bar{U}_{R_0}(t).
\end{align*}

Therefore,  for any $R_1, R_2>R_0$,
\begin{equation}
h_{R_1}(x,t)=h_{R_2}(x,t)=h_{R_0}(x,t).
\end{equation}
We complete the proof of Lemma \ref{heh}.
\End
Based on a similar argument as the one in the proof of Lemma \ref{heh}, we have the following lemma.
\begin{Lemma}[Uniform Invariance]\label{hehuu} For any $(x,t)\in \Sigma_k$ and for any $R_1, R_2>R'_0$, there holds
\begin{equation}
h_{R_1}(x,t)=h_{R_2}(x,t)=h_{R'_0}(x,t),
\end{equation}where $R'_0$ satisfies that  for any $R>R'_0$,
\begin{equation}\label{R00pie}
(\gm_R(s),h_R(\gm_R(s),s),\dot{\gm}_R(s))\in \mathcal{K}'\subset \mathcal{D}_{R'_0},
\end{equation}where $\gm_R:[0,t]\rightarrow M$ is a calibrated curve of $h_R$ satisfying $\gm_R(0)=x_0$ and $\gm_R(t)=x$..
\end{Lemma}

Next, we introduce the action function $h_{x_0,u_0}(x,t)$ for the original Lagrangian under (L1)-(L4) by a limiting passage. More precisely,
\begin{equation*}
h_{x_0,u_0}(x,t):=\lim_{R\rightarrow\infty}h_R(x,t).
\end{equation*}
By Lemma \ref{heh}, $h_{x_0,u_0}(x,t)$ is well-defined. By Lemma \ref{hehuu},  there exists $R'_0>0$ such that for any  $R>R'_0>0$, we have
\begin{equation*}
h_{x_0,u_0}(x,t)|_{\Sigma_k}=h_{R}(x,t).
\end{equation*}

\begin{Lemma}[Continuity]\label{continu}
$h_{x_0,u_0}(x,t)$ is continuous with respect to $(x,t)$ on $\Sigma_k$. In particular,  for $(\gm(t),t)\in \Sigma_k$, we have
\begin{equation}
\lim_{t\rightarrow 0^+} h_{x_0,u_0}(\gm(t),t)=u_0.
\end{equation}
\end{Lemma}
\Proof
Owing to Lemma \ref{exist} and Lemma \ref{unique}, $h_{x_0,u_0}(x,t)$ is continuous with respect to $(x,t)$ on $M\times (0,T]$. It follows that for any $\delta>0$, $h_{x_0,u_0}(x,t)$ is continuous with respect to $(x,t)$ on $M\times [\delta,T]$.
We denote
\[\Sigma_k^\epsilon:=\{(x,t)\in M\times[0,\epsilon]\ |\ d(x,x_0)\leq kt\},\]
where $d(\cdot,\cdot)$ denotes the distance between $x$ and $x_0$ induced by the Riemannian metric on $M$.

Based on Lemma \ref{xr} (ii), it follows from  a similar argument as Step (a) in Lemma \ref{chc} that for a given $k>0$, there exists $\epsilon>0$ such that for any  $(t,x)\in \Sigma_k^\epsilon$,
\begin{equation}
S(x,t)=h_{x_0,u_0}(x,t),
\end{equation}
where $S(x,t)$ denotes a $C^1$ classical solution of the Hamilton-Jacobi equation:
\begin{equation}
\partial_t S(x,t)+H(x,S(x,t),\partial_x S(x,t))=0.
\end{equation}
Hence, $h_{x_0,u_0}(x,t)$ is continuous with respect to $(t,x)\in \Sigma_k^\epsilon$. In order to prove Lemma \ref{continu}, it suffices to take $\delta<\epsilon$.
\End

\subsection{Proofs of Theorems}
In this step, we will prove Theorem \ref{two}, \ref{two2} and  \ref{two1} under the assumptions (L1)-(L4).

\noindent\textbf{Proof of Theorem \ref{two}:}
We will prove that
\begin{equation}\label{hxuxt}
h_{x_0,u_0}(x,t)=u_0+\inf_{\substack{\gm(t)=x \\  \gm(0)=x_0} }\int_0^tL(\gm(\tau),h_{x_0,u_0}(\gm(\tau),\tau),\dot{\gm}(\tau))d\tau,
\end{equation}where the infimum is taken among  the continuous and piecewise $C^1$ curves $\gm:[0,t]\rightarrow M$. Moreover, the infimum of (\ref{hxuxt}) is attained at a $C^1$ curve denoted by $\bar{\gm}$. Let
\[x(t):=\bar{\gm}(t),\quad u(t):=h_{x_0,u_0}(\bar{\gm}(t),t),\quad p(t)=\bar{\mathcal{L}}^{-1}(\bar{\gm}(t),h_{x_0,u_0}(\bar{\gm}(t),t),\dot{\bar{\gm}}(t)),\]
then $(x(t),u(t),p(t))$ satisfies the equation (\ref{hjech}).

According to Lemma \ref{inter}, there holds for each absolutely continuous curves $\gm:[0,t]\rightarrow M$, we have
\begin{equation}
h_R(x,t)\leq u_0+\int_0^tL_R(\gm(\tau),h_R(\gm(\tau),\tau),\dot{\gm}(\tau))d\tau.
\end{equation}
In particular, there exists a calibrated curve of $h_R$ denoted by $\gm_R: [0,t]\rightarrow M$ satisfying $\gm_R(0)=x_0$ and $\gm_R(t)=x$ such that
\begin{equation}\label{limitf}
h_R(x,t)= u_0+\int_0^tL_R(\gm_R(\tau),h_R(\gm_R(\tau),\tau),\dot{\gm}_R(\tau))d\tau.
\end{equation}

\vspace{1em}

\noindent\textbf{Claim A:}
 For each continuous and piecewise $C^1$ curves $\gm:[0,t]\rightarrow M$, there holds
\begin{equation}\label{hxuleq}
h_{x_0,u_0}(x,t)\leq u_0+\int_0^tL(\gm(\tau),h_{x_0,u_0}(\gm(\tau),\tau),\dot{\gm}(\tau))d\tau.
\end{equation}

\noindent\textbf{Proof of Claim A:}
By contradiction, we assume that there exists a continuous and piecewise $C^1$ curve $\bar{\gm}:[0,t]\rightarrow M$ such that
\begin{equation}\label{contradic}
h_{x_0,u_0}(x,t)> u_0+\int_0^tL(\bar{\gm}(\tau),h_{x_0,u_0}(\bar{\gm}(\tau),\tau),\dot{\bar{\gm}}(\tau))d\tau.
\end{equation}
Since $\bar{\gm}$ is a continuous and piecewise $C^1$ curve, then each component $\bar{\gm}$ has a bounded derivative which is continuous everywhere in $[0,t]$ except at a finite number of points at which left and right sided derivatives exist. Hence, it is easy to see that for any $s\in [0,t]$
\begin{equation}\label{hgam}
 |\dot{\bar{\gm}}(s\pm)|\leq C,
\end{equation}
where  $\dot{\bar{\gm}}(s\pm)$ denotes the left and right sided derivatives of $\bar{\gm}$ at $s$. Moreover, it follows from Lemma \ref{continu} that
\begin{equation}\label{hgam1}
|h_{x_0,u_0}(\bar{\gm}(s),s)|\leq C,
\end{equation}

Taking $R>\max\{C,R_0,R'_0\}$ where $R_0,R'_0$ are determined by (\ref{R00}) and (\ref{R00pie}) respectively, it follows from the construction of $L_R$ (see (\ref{llr})) that
\begin{align*}
h_{x_0,u_0}(x,t)=h_{R}(x,t)&\leq u_0+\int_0^tL_R(\bar{\gm}(\tau),h_R(\bar{\gm}(\tau),\tau),\dot{\bar{\gm}}(\tau))d\tau,\\
&=u_0+\int_0^tL(\bar{\gm}(\tau),h_{x_0,u_0}(\bar{\gm}(\tau),\tau),\dot{\bar{\gm}}(\tau))d\tau,
\end{align*}
which is in contradiction with (\ref{contradic}). Then (\ref{hxuleq}) is verified.\End

\noindent\textbf{Claim B:}
The infimum of (\ref{hxuleq}) can be attained at   a $C^1$ curve.

\noindent\textbf{Proof of Claim B:}
 By Lemma \ref{xr},  there exists a compact set $\mathcal{K}$  independent of $R$ such that for any $R$ and $s\in [0,t]$, there holds
\[\Phi_s^R(x_0,u_0,\dot{\gm}_R(0))=(\gm_R(s),h_R(\gm_R(s),s),\dot{\gm}_R(s))\in \mathcal{K},\]
where $\Phi_s^R$ denotes the contact flow generated by $L_R$. Hence, one can find a sequence $R_n$ such that as $n\rightarrow\infty$,
\begin{equation}
(x_0,u_0,\dot{\gm}_R(0))\rightarrow (x_0,u_0,v_0).
\end{equation}
Let $(\gm_\infty(s),u_\infty(s),\dot{\gm}_\infty(s)):=\Phi_s(x_0,u_0,v_0)$, where $\Phi_s$ denotes the contact flow generated by $L$. Based on the continuous dependence of solutions of ODEs on initial conditions, $\gm_\infty(t)=x$. Since  $\Phi_s(x_0,u_0,\dot{\gm}_{R_n}(0))$ converges uniformly on the compact interval $[0,t]$ to the map
\[s\mapsto (\gm_\infty(s),u_\infty(s),\dot{\gm}_\infty(s)),\]
Then one can pass the limit in (\ref{limitf}) to obtain
\begin{equation}
h_{x_0,u_0}(x,t)= u_0+\int_0^tL(\gm_\infty(\tau),h_{x_0,u_0}(\gm_\infty(\tau),\tau),\dot{\gm}_\infty(\tau))d\tau,
\end{equation}
which together with (\ref{hxuleq}) yields
\begin{equation}
h_{x_0,u_0}(x,t)=u_0+\inf_{\substack{\gm(t)=x \\  \gm(0)=x_0} }\int_0^tL(\gm(\tau),h_{x_0,u_0}(\gm(\tau),\tau),\dot{\gm}(\tau))d\tau,
\end{equation}where the infimum is taken among  the continuous and piecewise $C^1$ curves $\gm:[0,t]\rightarrow M$. In particular, the infimum  can be attained at  a $C^1$ curve.

\vspace{1em}

\noindent\textbf{Claim C:}
If the infimum of (\ref{hxuleq}) is attained  at  $\bar{\gm}$, then $(\bar{\gm}(\tau),u(\tau),p(\tau))$ satisfies the contact Hamiltonian equation (\ref{hjech}),
 where
\begin{equation}
u(\tau)=h(\bar{\gm}(\tau),\tau)\quad\text{and}\quad p(\tau)=\frac{\partial L}{\partial \dot{x}}(\bar{\gm}(\tau),h(\bar{\gm}(\tau),\tau),\dot{\bar{\gm}}(\tau)).
\end{equation}

\noindent\textbf{Proof of Claim C:}
By (\ref{hgam}) and (\ref{hgam1}), there holds
\begin{equation}
|h_{x_0,u_0}(\bar{\gm}(s),s)|\leq C, \quad |\dot{\bar{\gm}}(s\pm)|\leq C,
\end{equation}
where  $\dot{\bar{\gm}}(s\pm)$ denotes the left and right sided derivatives of $\bar{\gm}$ at $s$.
Take $R>\max\{C,R_0,R'_0\}$ where $R_0,R'_0$ are determined by (\ref{R00}) and (\ref{R00pie}) respectively.  Let $\gm_R:[0,t]\rightarrow M$ be a calibrated curve of $h_R$ satisfying $\gm_R(0)=x_0$ and $\gm_R(t)=x$. On one hand, it follows from the construction of $L_R$ (see (\ref{llr})) and Lemma \ref{xr} that
\begin{align*}
h_{x_0,u_0}(x,t)=h_{R}(x,t)&= u_0+\int_0^tL(\bar{\gm}(\tau),h_R(\bar{\gm}(\tau),\tau),\dot{\bar{\gm}}(\tau))d\tau,\\
&\leq u_0+\int_0^tL(\gm_R(\tau),h_R(\gm_R(\tau),\tau),\dot{\gm}_R(\tau))d\tau,\\
&=u_0+\int_0^tL_R(\gm_R(\tau),h_R(\gm_R(\tau),\tau),\dot{\gm}_R(\tau))d\tau,\\
&=u_0+\inf_{\substack{\gm(t)=x \\  \gm(0)=x_0} }\int_0^tL_R(\gm(\tau),h_R(\gm(\tau),\tau),\dot{\gm}(\tau))d\tau,\\
&=h_{R}(x,t).
\end{align*}
Hence, the third inequality is an equality. Meanwhile, we have
\[\int_0^tL(\bar{\gm}(\tau),h_R(\bar{\gm}(\tau),\tau),\dot{\bar{\gm}}(\tau))d\tau
=\int_0^tL_R(\bar{\gm}(\tau),h_R(\bar{\gm}(\tau),\tau),\dot{\bar{\gm}}(\tau))d\tau.\]
Moreover, we have
\begin{equation}\label{cl44}
\int_0^tL_R(\bar{\gm}(\tau),h_R(\bar{\gm}(\tau),\tau),\dot{\bar{\gm}}(\tau))d\tau=\inf_{\substack{\gm(t)=x \\  \gm(0)=x_0} }\int_0^tL_R(\gm(\tau),h_R(\gm(\tau),\tau),\dot{\gm}(\tau))d\tau,
\end{equation}
where $\gm\in C^{ac}([0,t],M)$. (\ref{cl44}) implies $\bar{\gm}$ is a calibrated curve of $h_R$. By virtue of Lemma \ref{chc}, $\bar{\gm}$ is $C^1$ and  $(\bar{\gm}(\tau),u(\tau),p(\tau))$ satisfies the contact Hamiltonian equation generated by $H_R$,
 where
\begin{equation}
u(\tau)=h_R(\bar{\gm}(\tau),\tau)\quad\text{and}\quad p(\tau)=\frac{\partial L_R}{\partial \dot{x}}(\bar{\gm}(\tau),h_R(\bar{\gm}(\tau),\tau),\dot{\bar{\gm}}(\tau)).
\end{equation}

Based on the construction of $L_R$ (see (\ref{llr})) and Lemma \ref{xr}, $(\bar{\gm}(\tau),u(\tau),p(\tau))$ also satisfies the contact Hamiltonian equation generated by $H$.
\End

So far, we complete the proof of Theorem \ref{two} under the assumptions (L1)-(L4).
\End
\begin{Remark}
In the proof of Theorem \ref{two}, in order to obtain (\ref{hgam}) more directly, we reduce the set of the absolutely continuous curves to the set of the continuous and piecewise $C^1$ curves for which the variational method can  still perform powerfully. Some technical difficulty would appear if more general set of the curves is considered.
\end{Remark}

\noindent\textbf{Proofs of  Theorem \ref{two2}:}
For any $R>R_0$ where $R_0$ is determined by (\ref{R00}), we have
\begin{equation*}
h_{x_0,u_0}(x,t)=h_{R}(x,t)=\bar{U}_R(t),
\end{equation*}where $(\bar{X}_{R}(t),\bar{U}_{R}(t),\bar{P}_{R}(t))$ denotes the solution of (\ref{hjech}) satisfying
\begin{equation}
\bar{U}_{R}(t)=\inf\left\{U_R(t):(X_R(s),U_R(s),P_R(s))\in \mathcal{S}_{R}\right\},
\end{equation}
where $\mathcal{S}_{R}$ denotes the set of solutions  $(X_{R}(s),U_{R}(s),P_{R}(s))$  of (\ref{hjech}) generated by $H_R$ satisfying $X_{R}(0)=x_0$, $X_{R}(t)=x$ and $U_{R}(0)=u_0$.
In order to prove Theorem \ref{two2}, it suffices to show
\begin{equation}
\bar{U}_{R}(t)=\inf\left\{U(t):(X(s),U(s),P(s))\in \mathcal{S}\right\},
\end{equation}
where $\mathcal{S}$ denotes the set of solutions of (\ref{hjech}) $(X(s),U(s),P(s))$ satisfying $X(0)=x_0$, $X(t)=x$ and $U(0)=u_0$.
Based on the construction of $L_R$ (see (\ref{llr})) and Lemma \ref{xr}, there holds
\[(\bar{X}_{R}(t),\bar{U}_{R}(t),\bar{P}_{R}(t))\in  \mathcal{S}_{R}\cap\mathcal{D}^*_R=\mathcal{S}\cap\mathcal{D}^*_R\subset \mathcal{S},\]
which implies
\begin{equation}\label{urttt}
\bar{U}_{R}(t)\geq\inf\left\{U(t):(X(s),U(s),P(s))\in \mathcal{S}\right\}.
\end{equation}

It remains to show that for each $(X(s),U(s),P(s))\in \mathcal{S}$,  $U(t)\geq \bar{U}_{R}(t)$, from which it yields that
\begin{equation}\label{urttt22}
\bar{U}_{R}(t)\leq\inf\left\{U(t):(X(s),U(s),P(s))\in \mathcal{S}\right\}.
\end{equation}
By contradiction, we assume that there exists $(\tilde{X}(t),\tilde{U}(t),\tilde{P}(t))\in \mathcal{S}$ such that
\begin{equation}\label{hhhxx}
\tilde{U}(t)<\bar{U}_{R}(t).
\end{equation}
Since $|\tilde{U}(s)|$ and $|\tilde{P}(s)|$ have a bound denoted by $C$ for any $s\in [0,t]$, then
\[(\tilde{X}(s),\tilde{U}(s),\tilde{P}(s))\in \mathcal{D}^*_C.\]
 Take $R>\max\{C,R_0,R'_0\}$ where $R_0,R'_0$ are determined by (\ref{R00}) and (\ref{R00pie}) respectively, we have
 \begin{align*}
\tilde{U}(t)&\geq \inf\left\{U(t):(X(s),U(s),P(s))\in \mathcal{S}\cap\mathcal{D}^*_C\right\},\\
 &=\inf\left\{U_R(t):(X_R(s),U_R(s),P_R(s))\in \mathcal{S}_R\cap\mathcal{D}^*_C\right\},\\
 &\geq\inf\left\{U_R(t):(X_R(s),U_R(s),P_R(s))\in \mathcal{S}_R\right\},\\
 &=\bar{U}_R(t),
 \end{align*}
which is in contradiction with (\ref{hhhxx}).
Therefore, we have
\begin{equation*}
h_{x_0,u_0}(x,t)=\bar{U}_{R}(t)=\inf\left\{U(t):(X(s),U(s),P(s))\in \mathcal{S}\right\},
\end{equation*}
where $\mathcal{S}$ denotes the set of solutions of (\ref{hjech}) $(X(s),U(s),P(s))$ satisfying $X(0)=x_0$, $X(t)=x$ and $U(0)=u_0$.

So far, we complete the proofs of  Theorem \ref{two2} under the assumptions (L1)-(L4).
\End

\noindent\textbf{Proofs of  Theorem \ref{two1}:}
Theorem \ref{two1} asserts that
for each solution $(X(t),U(t),P(t))$ of  (\ref{hjech}), there exists $\epsilon>0$ such that for any $t\in [0,\epsilon]$,
\[U(t)=h_{x_0,u_0}(X(t),t),\]
where $(X(t),U(t),P(t)$ denotes the solution of the equation (\ref{hjech}) satisfying $X(0)=x_0$, $U(0)=u_0$ and $P(0)=p_0$. We omit the subscript $x_0,u_0$ of $h_{x_0,u_0}$ for simplicity.

Let $\bar{\gm}:[0,t]\rightarrow M$ be the calibrated curve of $h$ with $\bar{\gm}(0)=x_0$ and $\bar{\gm}(t)=X(t)$. By virtue of Theorem \ref{two}, $(\bar{\gm}(s),h(\bar{\gm}(s),s),\dot{\bar{\gm}}(s))$ coincidences with $(\bar{X}(s),\bar{U}(s),\dot{\bar{X}}(s))$ for $s\in [0,t]$, where $(\bar{X}(s),\bar{U}(s),\dot{\bar{X}}(s))$ denotes the dual variables of $(\bar{X}(s),\bar{U}(s),\bar{P}(s))$ via the Legendre transformation.

Since $\bar{U}(t)=h(\bar{X}(t),t)$, it follows from Theorem \ref{two2} that
\begin{equation}\label{84yy1}
\bar{U}(t)\leq U(t)\leq C.
\end{equation}
By a similar argument as Lemma \ref{aprioriset}, it yields that for any $s\in [0,t]$, $|\bar{U}(s)|<C$, where $C$ is a constant independent of $s$.

Since $(X(t),U(t),P(t))$ is the solution of the equation (\ref{hjech}) satisfying $X(0)=x_0$, $U(0)=u_0$ and $P(0)=p_0$, then $|P(s)|\leq C$ for any $s\in [0,t]$, then $|\dot{X}(s)|\leq C$. Hence, we have
\begin{equation}\label{xrscbyy}
|U(t)-u_0|\leq Ct.
\end{equation}
According to  (\ref{84yy1}), there holds
\begin{align*}
\int_{0}^{t}L(\bar{X}(s),\bar{U}(s),\dot{\bar{X}}(s))ds\leq U(t)-u_0\leq Ct.
\end{align*}
It follows from (L2) that
\[\int_{0}^{t}|\dot{\bar{X}}(s)|-Bds\leq Ct.\]
By the mean value theorem for integration, that there exists $s_0\in [0,t]$ such that
\begin{equation}\label{xrscb99}
|\dot{\bar{X}}(s_0)|\leq B+C.
\end{equation}
Based on the continuous dependence of solutions of ODEs on initial conditions and (L3), it follows that for any $s\in [0,t]$,
\[(\bar{X}(s),\bar{U}(s),\dot{\bar{X}}(s))\in \mathcal{K}.\]

 Let $X(t,v)$ be the first argument of  $\Phi_{t}(x_0,u_0,v)$, where $\Phi$ denotes the flow generated by $L$. It follows from the implicit function theorem that there exists  $\epsilon>0$ small enough such that for $t\in [0,\epsilon]$,  $X(t,\cdot)$ is a diffeomorphism onto its image. Hence, if $\bar{X}(t)=X(t)$ for $t\leq \epsilon$,
then
\[\dot{\bar{X}}(0)=\dot{X}(0),\]
which together with $\bar{X}(0)=X(0)=x_0$ and $\bar{U}(0)=U(0)=u_0$ implies that for any $s\in [0,t]$,
\[U(s)=h(X(s),s).\]

So far, we have completed the proof of Theorem \ref{two1}.\End

 \vspace{2ex}
\noindent\textbf{Acknowledgement}
 % The authors sincerely
%thanks the referees for their careful reading of the manuscript and
%invaluable comments which were very helpful in improving this paper.
%The authors  would like to thank
%Prof. C.-Q. Cheng and Dr. L. Jin for many helpful discussions.
L. Wang is partially under the support of National Natural Science Foundation of China (Grant No. 11401107). J. Yan is partially under the support of National Natural Science Foundation of China (Grant No. 11171071,  11325103) and
National Basic Research Program of China (Grant No. 2013CB834100).

\vspace{2em}

{\sc Lin Wang}

{\sc School of Mathematical Sciences, Fudan University,
Shanghai 200433,
China.}

 {\it E-mail address:} \texttt{linwang.math@gmail.com}

\vspace{1em}

{\sc Jun Yan}

{\sc School of Mathematical Sciences, Fudan University,
Shanghai 200433,
China.}

 {\it E-mail address:} \texttt{yanjun@fudan.edu.cn}

\end{document}